\documentclass[11pt,a4paper]{article}

\usepackage[utf8]{inputenc}
\usepackage[T1]{fontenc}
\usepackage[margin=1in]{geometry}
\usepackage{lmodern}
\usepackage{microtype}
\usepackage{graphicx}
\usepackage{amsmath,amssymb,amsthm}
\usepackage{mathtools}
\usepackage{authblk}
\usepackage{placeins}
\usepackage{url}
\usepackage{appendix}
\usepackage{natbib}
\usepackage{xcolor}
\usepackage{pgfplots}
\usepackage{tikz}
\usetikzlibrary{arrows.meta,positioning,decorations.pathreplacing,decorations.markings}
\pgfplotsset{compat=1.18}

\newtheorem{theorem}{Theorem}
\newtheorem{lemma}{Lemma}
\newtheorem{proposition}{Proposition}
\newtheorem{remark}{Remark}

\newcommand{\K}{\mathcal{K}}

\newcommand{\pS}{p_{\mathrm{S}}}
\newcommand{\pI}{p_{\mathrm{I}}}
\newcommand{\pR}{p_{\mathrm{R}}}

\newcommand{\dpS}{\dot{p}_{\mathrm{S}}}

\newcommand{\dpR}{\dot{p}_{\mathrm{R}}}
\newcommand{\Sinf}{S_{\infty}}
\newcommand{\Rinf}{R_{\infty}}
\newcommand{\Rz}{\mathcal{R}_0}
\newcommand{\diff}{\mathrm{d}}
\newcommand{\ind}{i_0}
\newcommand{\Pop}{P_L}

\title{A Low-Dimensional SIR Model \\ for Epidemic Propagation in a Network}
\author[1]{M. Soledad Aronna}
\author[2,3]{Mariana Bergonzi}
\author[2,3]{Ernesto Kofman}
\affil[1]{School of Applied Mathematics, FGV EMAp, Brazil}
\affil[2]{French-Argentine International Center for Information and System Sciences (CIFASIS), CONICET, Argentina}
\affil[3]{FCEIA, Universidad Nacional de Rosario, Argentina}
\date{}

\begin{document}

\maketitle

\begin{abstract}
We propose an extension of the classical susceptible--infectious--recovered (SIR) model that includes two additional compartments designed to capture the dynamics of disease spread across multiple interconnected cities or populated regions, while avoiding the high dimensionality and large parameter sets typical of network-based or agent-based approaches. Instead of explicitly modeling individual locations or mobility networks, we introduce aggregate variables that describe whether the epidemic has not yet reached, is currently active in, or has already passed through different regions of the spatial domain. This formulation allows the model to reproduce key qualitative features observed in aggregated incidence data, such as prolonged plateaus and multiple infection waves arising from asynchronous local outbreaks. The resulting system consists of five ordinary differential equations with a small number of interpretable parameters, providing a tractable framework for analytical investigation and numerical simulation. Our approach offers a parsimonious alternative for studying spatially structured epidemic dynamics when only aggregated data are available or when model simplicity is essential.
\end{abstract}


\section{Introduction}

The aim of this work is to extend the classical SIR (susceptible--infectious--recovered) model by including new compartments that capture the level of \emph{spatial propagation} of the epidemic. The main objective is to find a system that can imitate the dynamics observed when a disease spreads across different interconnected cities or populated regions, while---unlike network-based models (see references below)---involving a relatively small number of parameters.

Models for the spatial propagation of diseases have existed for quite a long time. Early work by \cite{may1984spatial} from the 1980s considered the impact of spatial heterogeneity on the vaccination requirement for the eradication of a canonical infectious disease. That study built upon earlier contributions by \cite{hethcote1978immunization} and \cite{post1983endemic}, which highlighted the importance of heterogeneity in epidemiological modeling. Further references and historical context on foundational articles can be found in \cite{rhodes1996persistence} and \cite{bolker1995space}.

The works of \cite{sattenspiel1995structured} and \cite{arino2003multi} explored how regional mobility and migration can be incorporated into epidemic models. In particular, \cite{arino2003multi} emphasized the use of space-discrete models and derived inequalities for the basic reproduction number, $\mathcal{R}_0$. Considering the structure of urban connectivity networks, \cite{colizza2007reaction} and \cite{colizza2008epidemic} developed a framework to calculate global thresholds for epidemic invasion. Complementary to these works, \cite{pastor2001epidemic} and \cite{pastor2002epidemic} made seminal contributions to the study of epidemic dynamics on complex networks of individuals. Recently, \cite{aronna2024optimal} conducted a study addressing how intercity commuting patterns within metropolitan regions influence the spread of infectious diseases.

Agent-based simulation models provide an alternative to equation-based systems and, in some cases, can represent richer scenarios than the latter \citep{hunter2018comparison}. In an agent-based framework, each individual agent can be assigned distinct attributes and decision rules, enabling the model to capture interactions and behaviors at the individual level \citep{marshall2015formalizing,hunter2017taxonomy}. In contrast, the approach taken in this article differs substantially from agent-based modeling: the defining feature of our model is aggregation rather than individual-level representation.

In this article, we propose an ODE model that extends the classical three-state SIR framework by adding two state variables, $p_S$ and $p_R$, that represent spatial propagation. We call it the \emph{pSIR model}. Loosely speaking, the total population is assumed to live in several interconnected districts and is partitioned according to whether the epidemic is active in each location, has not yet arrived ($p_S$), or has already passed ($p_R$). We then define the dynamics of the resulting five-state model so that it approximates the aggregated behavior of a metapopulation model on a network of arbitrary topology. The goal is to capture the essential features of spatial diffusion across a network of cities or regions, where the epidemic begins in one location and subsequently spreads to others, generating prolonged plateaus or multiple infectious waves.

The article is organized as follows. In Section~\ref{SecModel}, we introduce the structure of the model and deduce the expressions for its variables and parameters. In Section~\ref{sec:properties}, we discuss the fundamental mathematical properties of the model, including well-posedness, asymptotic behavior, final epidemic size relations, and the calculation of the basic reproduction number $\Rz$. Section~\ref{SecExamples} presents numerical simulations that corroborate the consistency between the results obtained using the pSIR model and those provided by high-dimensional metapopulation models, as well as a fit to COVID-19 data from Argentina. Finally, Section~\ref{sec:conclusions} concludes the manuscript and proposes future lines of research.

\section{Model Formulation} \label{SecModel}

In this section, we introduce a novel epidemic model. We first propose the structure of the model and then deduce the dynamics of each variable in such a way that it captures the aggregated behavior of a metapopulation model.

\subsection{Aggregated pSIR Model Structure}

We consider an epidemic taking place in a population of size $N$ that is divided into a certain number of districts. Individuals remain most of the time in their own district, but they sometimes move to other districts where they can get infected or, if they are already infected, they can infect other individuals.

We first divide the population into three main groups $\pS$, $\pI$, and $\pR$, defined as follows: $\pS(t)$ is the population living in districts that have not yet been reached by the epidemic at time $t$; $\pI(t)$ is the population living in districts where the epidemic is active at time $t$; and $\pR(t)$ is the population living in districts where the epidemic has already finished. Thus,
\begin{equation} \label{eq:N}
\pS(t) + \pI(t) + \pR(t) = N.
\end{equation}
The population $\pS(t)$ consists exclusively of susceptible individuals, whereas $\pI(t)$ comprises susceptible, infected, and recovered individuals, denoted by $S(t)$, $I(t)$, and $R(t)$, respectively, so that
\begin{equation} \label{eq:pIdef}
S(t) + I(t) + R(t) = \pI(t).
\end{equation}
Finally, the population $\pR(t)$ is composed of recovered individuals as well as susceptible ones, since not every individual is infected before the local outbreak ends.

\begin{figure}[htb]
\centering
\includegraphics[width=5cm]{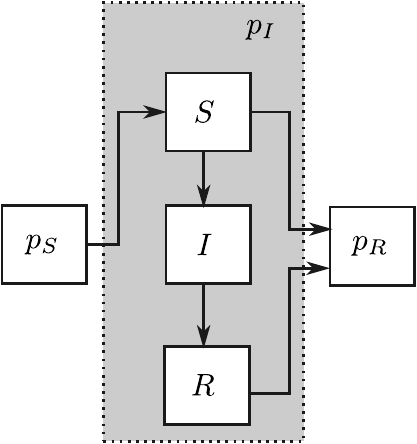}
\caption{Aggregated pSIR model scheme.}
\label{fig:psir}
\end{figure}

The three groups are coupled through two flows. As the epidemic reaches new districts, susceptible individuals leave $\pS(t)$ and enter $\pI(t)$ at rate $-\dpS(t)$. Likewise, as districts recover---that is, as the wave moves away from them---individuals leave $\pI(t)$ and enter $\pR(t)$ at rate $\dpR(t)$. To split this second flow, we let $\Sinf$ denote the fraction of individuals in a district who remain susceptible once the epidemic has finished there, and $\Rinf \coloneq 1 - \Sinf$ the fraction who were infected. Susceptible and recovered individuals then move toward $\pR(t)$ at rates $\Sinf\,\dpR(t)$ and $\Rinf\,\dpR(t)$, respectively. Figure~\ref{fig:psir} shows the compartmental diagram of the model.

We assume that within the population $\pI(t)$, the dynamics is governed by a modified SIR model that accounts for the arrival of susceptible individuals from $\pS$ at rate $-\dpS(t)$ and the departure of susceptible and recovered individuals toward $\pR$. Thus, the dynamics for the population $\pI(t)$ can be written as
\begin{equation} \label{eq:SIRp}
\begin{split}
\dot{S}(t) &= -\dpS(t) - \beta \frac{S(t) I(t)}{\pI(t)} - \Sinf \dpR(t), \\
\dot{I}(t) &= \beta \frac{S(t) I(t)}{\pI(t)} - \gamma I(t), \\
\dot{R}(t) &= \gamma I(t) - \Rinf \dpR(t).
\end{split}
\end{equation}
To complete the model, we need to establish expressions for $\dpS(t)$ and $\dpR(t)$ representing the advance of the epidemic through the different districts and their subsequent recovery.

\subsection{Derivation of $\dot{p}_S$}

We want the pSIR model of Eq.~\eqref{eq:SIRp} to approximate the aggregated dynamics of a metapopulation model consisting of $D$ districts, whose individual epidemic dynamics, for every $i = 1, \dots, D$, are described by
\begin{equation}\label{eq:SIRi}
\begin{split}
\dot{S}_i(t) &= -\beta \frac{S_i(t) I_i(t) (1 - c_i p)}{N_d} - \beta \frac{S_i(t)}{N_d} \sum_{j \in C_i} p I_j(t), \\
\dot{I}_i(t) &= \beta \frac{S_i(t) I_i(t) (1 - c_i p)}{N_d} + \beta \frac{S_i(t)}{N_d} \sum_{j \in C_i} p I_j(t) - \gamma I_i(t), \\
\dot{R}_i(t) &= \gamma I_i(t),
\end{split}
\end{equation}
where $C_i$ denotes the set of neighboring districts of locality $i$, $c_i \coloneq |C_i|$ is the corresponding number of neighbors, $p$ represents the daily time proportion that a resident of a region spends in neighboring regions, and $N_d \coloneq N/D$ is the population of each district.

\begin{figure}[htb]
\centering
\begin{tikzpicture}[
  scale=0.85,
  d/.style={circle, draw, fill=white, inner sep=2pt},
  ini/.style={circle, draw, fill=red!60, inner sep=2.6pt},
  every edge/.style={draw, gray!80}
]
  \node[ini, label=below:{$i_0$}] (a) at (0,0) {};
  \node[d] (b1) at (1.8, 1.1) {};
  \node[d] (b2) at (1.8, 0.0) {};
  \node[d] (b3) at (1.8,-1.1) {};

  \node[d] (c1) at (3.6, 1.5) {};
  \node[d] (c2) at (3.6, 0.5) {};
  \node[d] (c3) at (3.6,-0.5) {};
  \node[d] (c4) at (3.6,-1.5) {};

  \node[d] (d1) at (5.4, 1.0) {};
  \node[d] (d2) at (5.4, 0.0) {};
  \node[d] (d3) at (5.4,-1.0) {};

  \node[d] (e1) at (7.2, 0.55) {};
  \node[d] (e2) at (7.2,-0.55) {};

  \node[font=\small] at (0, 2.15) {$L_0$};
  \node[font=\small] at (1.8, 2.15) {$L_1$};
  \node[font=\small] at (3.6, 2.15) {$L_2$};
  \node[font=\small] at (5.4, 2.15) {$L_3$};
  \node[font=\small] at (7.2, 2.15) {$L_4$};

  \draw[gray!75] (a)--(b1) (a)--(b2) (a)--(b3);
  \draw[gray!75] (b1)--(c1) (b1)--(c2) (b2)--(c2) (b2)--(c3) (b3)--(c4);
  \draw[gray!75] (c1)--(d1) (c2)--(d1) (c3)--(d2) (c3)--(d3) (c4)--(d3);
  \draw[gray!75] (d1)--(e1) (d2)--(e1) (d2)--(e2) (d3)--(e2);

  \draw[gray!75] (b1)--(b2) (c2)--(c3) (c3)--(c4) (d1)--(d2);
\end{tikzpicture}
\caption{Example of a district network organized in layers: $L_0=\{i_0\}$ is the initial district, and $L_k$ contains the districts at distance $k$ from $i_0$. Edges join districts in consecutive layers or in the same layer.}
\label{FigLk}
\end{figure}
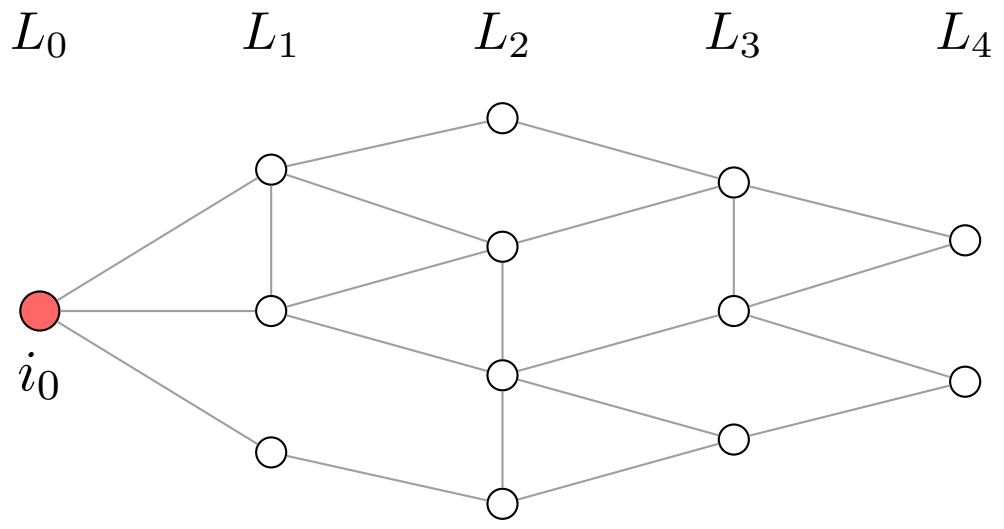

Assume that at time $t = 0$, the epidemic begins in a single district $\ind$, and let $L_0 \coloneq \{\ind\}$. We define the layer $L_1$ as the set of districts that are neighbors of $\ind$, and, in general, the layer $L_k$ as the set of districts at graph distance $k$ from $\ind$ (see Figure~\ref{FigLk}).

We assume that the epidemic propagates through consecutive layers such that the population in layer $L_k$ moves from $p_S$ to $p_I$ before the population in layer $L_{k+1}$ does so. The layer currently transitioning from $p_S$ to $p_I$ is called the \emph{boundary layer}. The size of the population in the boundary layer at a given time can be computed as a function $\Pop$ of the population remaining in $p_S$. We assume that $\Pop(p_S)$ is a known function characterized by the network topology; a generic expression for $\Pop$ is derived later in Section~\ref{SubsecpL}.

The epidemic is considered to propagate between consecutive layers at a constant mean rate. That is, every $\Delta t$ units of time, an entire layer leaves $p_S$, resulting in
\begin{equation*}
p_S(t + \Delta t) = p_S(t) - \Pop(p_S).
\end{equation*}
We can therefore approximate
\begin{equation} \label{eq:dotpS1}
\dot{p}_S(t) \approx -\frac{\Pop(p_S)}{\Delta t}.
\end{equation}
This approximation is only valid after the infected population $I$ reaches a certain threshold, since the time at which the first layer $L_1$ is reached depends on the number of initially infected individuals and is not necessarily equal to $\Delta t$. Thus, for the early stage of the epidemic, we next analyze the dynamics of the index district $\ind$ and its neighboring districts in the first layer $L_1$.

We require $\dot{I}(t) \approx \dot{I}_{\ind} + \sum_{i \in L_1} \dot{I}_i$ so that the pSIR model effectively approximates the aggregated behavior. From Eqs.~\eqref{eq:SIRp} and~\eqref{eq:SIRi}, it follows that
\begin{equation}\label{eq:startCond}
\begin{split}
\beta I \frac{S}{p_I} - \gamma I
&\approx \beta I_{\ind} \frac{S_{\ind}}{N_d} - \gamma I_{\ind} + \sum_{i \in L_1} \left(\beta I_i \frac{S_i}{N_d} - \gamma I_i\right) \\
&\approx \beta I_{\ind} \frac{S_{\ind}}{N_d} - \gamma I_{\ind} + |L_1| \left(\beta I_{i_1} \frac{S_{i_1}}{N_d} - \gamma I_{i_1}\right),
\end{split}
\end{equation}
where $i_1 \in L_1$, $|L_1|$ denotes the number of districts in layer $L_1$, and we assume that all districts in $L_1$ have identical evolution. We also consider that initially $I_{i_1} \approx 0$, so the cross-infection term $\sum_{j \in C_0} p I_j(t)$ is neglected.

At the beginning, while $I_{\ind} \gg \sum_{i \in L_1} I_i$ and $\dot{I}_{\ind} \gg \sum_{i \in L_1} \dot{I}_i$, condition~\eqref{eq:startCond} is satisfied with $p_I \approx N_d$, $I \approx I_{\ind}$, and $S \approx S_{\ind}$. This implies that we can select $p_I(0) = N_d$ (an entire district inside compartment $p_I$) and must keep $p_I$ almost unchanged while $I_{\ind} \gg \sum_{i \in L_1} I_i$.

After a certain time, the entire population of layer $L_1$ moves out of $\pS$, and $\pI$ reaches $(|L_1|+1)N_d$, with $S = S_{\ind} + \sum_{i \in L_1} S_i$ and $I = I_{\ind} + \sum_{i \in L_1} I_i$. Substituting these relations into~\eqref{eq:startCond}, we obtain
\begin{equation*}
\beta I \frac{S}{(|L_1|+1)N_d} - \gamma\bigl(I_{\ind} + |L_1|I_{i_1}\bigr)
= \beta I_{\ind}\frac{S_{\ind}}{N_d} + |L_1|\beta I_{i_1}\frac{S_{i_1}}{N_d} - \gamma\bigl(I_{\ind} + |L_1|I_{i_1}\bigr).
\end{equation*}
Simple algebraic manipulation of this equation leads to
\begin{equation*}
\bigl(I_{\ind} - I_{i_1}\bigr)\bigl(S_{i_1} - S_{\ind}\bigr) = 0.
\end{equation*}
Since $S_{i_1} \neq S_{\ind}$, we conclude that $I_{\ind} = I_{i_1}$. This means that when $I_{\ind} = I_{i_1}$ (a situation that can only occur after district $\ind$ passes its peak), we must have $p_I = (|L_1|+1)N_d$, implying that $p_I$ must begin increasing as $I(t)$ approaches the peak size $I_{\mathrm{peak}}$ of the initial district. This can be achieved by modifying Eq.~\eqref{eq:dotpS1} as follows:
\begin{equation} \label{eq:dpS}
\dot{p}_S(t) = -\frac{\Pop(p_S)}{\Delta t} \min\!\left(1, (I/I_{\mathrm{peak}})^2\right),
\end{equation}
where $I_{\mathrm{peak}}$ is given by the classical SIR peak formula \cite{hethcote2000mathematics}:
\begin{equation}\label{eq:Ipeak}
I_{\mathrm{peak}} = N_d \left(1 - \frac{1 + \log \Rz}{\Rz}\right).
\end{equation}

\subsubsection{Derivation of $\Delta t$}

We next estimate the elapsed time $\Delta t$ from the moment the epidemic reaches one layer until it reaches the next. To this end, we assume that the districts of the $k$-th layer are \emph{reached} at time $t_k$, defined as the time at which $I_i(t)$ reaches a threshold $I_{\mathrm{th}}$ for $i \in L_k$. By taking $I_{\mathrm{th}} \ll I_{\mathrm{peak}}$, we can approximate $S_i \approx N_d$, so that Eq.~\eqref{eq:SIRi} simplifies to
\begin{equation} \label{eq:Iapprox}
\dot{I}_i(t) \approx \beta I_i(t) + \beta \sum_{j \in C_i} p I_j(t) - \gamma I_i(t).
\end{equation}
Provided that $p \ll 1$ and $I_{\mathrm{th}}$ is large enough that $p I_j(t_k) \ll I_{\mathrm{th}}$, Eq.~\eqref{eq:Iapprox} reduces to $\dot{I}_i(t) \approx (\beta - \gamma)I_i(t)$, yielding
\begin{equation*}
I_i(t) \approx e^{\lambda (t - t_k)} I_{\mathrm{th}}
\end{equation*}
with $\lambda = \beta - \gamma$ for $t > t_k$. This exponential expression also holds for the trajectories of a district $j$ in the $(k-1)$-th layer:
\begin{equation} \label{eq:Ijexp}
I_j(t) \approx e^{\lambda (t - t_{k-1})} I_{\mathrm{th}},
\end{equation}
where $t_{k-1} = t_k - \Delta t$ is the time at which the $(k-1)$-th layer reaches the threshold $I_{\mathrm{th}}$. Substituting this into Eq.~\eqref{eq:Iapprox} gives
\begin{equation*}
\dot{I}_i(t) \approx \lambda I_i(t) + \beta p |M_i| e^{\lambda (t - t_{k-1})} I_{\mathrm{th}},
\end{equation*}
where $|M_i| = |C_i \cap L_{k-1}|$ is the number of neighbors of district $i$ belonging to the previous layer. Solving this linear differential equation with $I_i(t_{k-1}) \approx 0$ (since at $t_{k-1}$ the $k$-th layer had not yet been reached) yields
\begin{equation} \label{eq:Itexp}
I_i(t) \approx \beta p |M_i| (t - t_{k-1}) e^{\lambda(t - t_{k-1})} I_{\mathrm{th}}.
\end{equation}
Since $I_i(t)$ is delayed by $\Delta t$ with respect to $I_j(t)$, they must satisfy $I_i(t_{k+1}) = I_j(t_k)$. Substituting Eqs.~\eqref{eq:Ijexp} and~\eqref{eq:Itexp} gives
\begin{equation*}
\beta p |M_i| (t_{k+1} - t_{k-1}) e^{\lambda(t_{k+1} - t_{k-1})} I_{\mathrm{th}} = e^{\lambda (t_k - t_{k-1})} I_{\mathrm{th}},
\end{equation*}
which implies
\begin{equation*}
\beta p |M_i| 2 \Delta t e^{\lambda \Delta t} = 1.
\end{equation*}
The solution for $\Delta t$ can be expressed in terms of the Lambert $\mathrm{W}$ function as
\begin{equation*}
\Delta t = \frac{1}{\lambda} \mathrm{W}\!\left(\frac{\lambda}{2 m p \beta}\right) = \frac{1}{\beta - \gamma} \mathrm{W}\!\left(\frac{\beta - \gamma}{2 m p \beta}\right),
\end{equation*}
where $m$ is the average value of $|M_i|$. If districts have, on average, $c$ connections, we can take $m = c/2$, assuming that half of the connections lead to the previous layer and the other half to the next layer. Then, the final expression for $\Delta t$ is
\begin{equation} \label{eq:dt}
\Delta t = \frac{1}{\beta - \gamma} \mathrm{W}\!\left(\frac{\beta - \gamma}{c p \beta}\right),
\end{equation}
which depends on the parameters $\beta$, $\gamma$, the mobility $p$, and the connectivity $c$.

\subsubsection{Computation of $\Pop(p_S)$} \label{SubsecpL}

In this subsection, we deduce an expression for the function $\Pop$, which gives the number of individuals in the boundary layer. Let $P(k)$ denote the size of the population residing in layer $L_k$ for $k = 0, 1, \dots, K$, with $P(k) = 0$ for $k > K$. The cumulative population up to layer $L_k$ is given by
\begin{equation*}
P_c(k) = \sum_{j=0}^k P(j),
\end{equation*}
and we set $P_c(-1) = 0$.

\begin{figure}[ht]
\centering
\begin{tikzpicture}[x=1.4cm,y=0.55cm,>=stealth,scale=0.8,transform shape]
\def\NpS{6.3}
\def\a{4.5}
\def\b{7}
\def\N{10.5}
\fill[red!8] (0.0,\a) rectangle (6.2,\b);
\draw[red!50,densely dotted] (0,\a) -- (6.2,\a) (0,\b) -- (6.2,\b);
\draw[->] (0,0) -- (6.4,0) node[below right] {layer $L_j$};
\draw[->] (0,0) -- (0,11.6) node[above] {cumulative population};
\draw[very thick,blue!70!black] (0,0) -- (0,1) -- (1,1) -- (1,2.5) -- (2,2.5) -- (2,4.5) -- (3,4.5) -- (3,7) -- (4,7) -- (4,9) -- (5,9) -- (5,10.5) -- (6,10.5);
\draw[line width=2.2pt,red!75!black] (3,\a) -- (3,\b);
\foreach \x/\y in {0/1,1/2.5,2/4.5,3/7,4/9,5/10.5}
  \fill[blue!70!black] (\x,\y) circle (2.2pt);
\foreach \x/\y in {0/0.5,1/1.75,2/3.5,3/5.75,4/8,5/9.75}
  \filldraw[fill=white,draw=blue!70!black,thick] (\x,\y) circle (2.2pt);
\draw[thick,dashed] (0,\NpS) -- (6.2,\NpS) node[right] {$N-p_S$};
\fill[red!75!black] (3,\NpS) circle (2.6pt);
\foreach \x/\lab in {0/0,1/1,2/k-1,3/k,4/k+1,5/K}
  \draw (\x,0) -- (\x,-0.2) node[below] {$\lab$};
\foreach \y/\lab in {\a/P_c(k-1),5.75/m_k,\b/P_c(k),\N/N=P_c(K)}
  \draw (0,\y) -- (-0.08,\y) node[left,font=\small] {$\lab$};
\draw[densely dotted] (0,5.75) -- (3,5.75);
\draw[densely dotted] (0,\N) -- (5,\N);
\draw[decorate,decoration={brace,mirror,amplitude=5pt}] (7.3,0) -- (7.3,\NpS-0.05) node[midway,right=6pt,font=\small] {$p_I+p_R$};
\draw[decorate,decoration={brace,mirror,amplitude=5pt}] (7.3,\NpS+0.05) -- (7.3,\N) node[midway,right=6pt,font=\small] {$p_S$};
\draw[<-,thick,red!75!black] (3.08,5.0) -- (4.1,2.6) node[below right,align=left,font=\small,text=red!75!black,inner sep=1pt] {$N-p_S\in[P_c(k-1),P_c(k)]$\\ $\Rightarrow$ $L_k$ is the boundary layer};
\begin{scope}[shift={(0.25,10.9)},font=\small]
\fill[blue!70!black] (0,0) circle (2.2pt);
\node[right] at (0.08,0) {$P_c(j)$};
\filldraw[fill=white,draw=blue!70!black,thick] (1.1,0) circle (2.2pt);
\node[right] at (1.18,0) {midpoint $m_j$};
\end{scope}
\end{tikzpicture}
\caption{Cumulative population $P_c(j)$ as a step function of the layer index. Each step spans the range $[P_c(j-1),P_c(j)]$, of height $P(j)$, with midpoint $m_j$ (hollow dots). The value $N-p_S$ (dashed line) crosses the step of layer $L_k$, which is therefore the boundary layer; its midpoint $m_k$ yields $\Pop(p_S=N-m_k)=P(k)$.}
\label{fig:boundary_layer}
\end{figure}
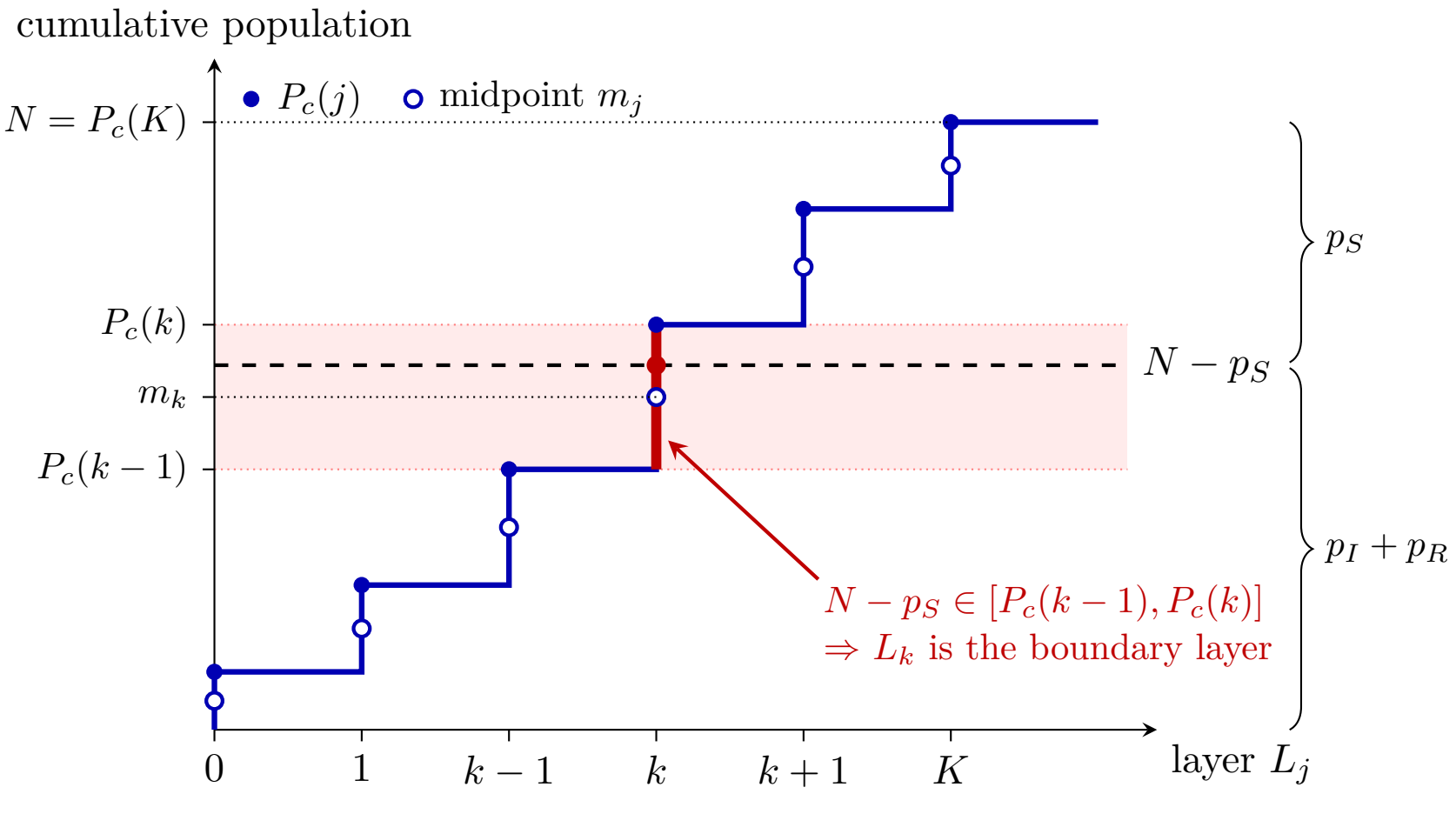

Notice that $L_k$ is the boundary layer when $N - p_S$ lies within the accumulated range $[P_c(k-1), P_c(k)]$ (see Figure~\ref{fig:boundary_layer}). Thus, at the midpoint $m_k = (P_c(k-1) + P_c(k))/2$ of this interval, we set
\begin{equation} \label{eq:Popdisc}
\Pop(\pS = N - m_k) = P(k),
\end{equation}
which defines $\Pop$ at these discrete points.

\subsubsection*{Interpolation Formula for $\Pop(p_S)$}

\begin{figure}[ht]
\centering
\begin{tikzpicture}[x=1cm,y=1cm,>=stealth,scale=0.85,transform shape]
\foreach \l/\r/\h/\lab in {9.5/10.5/1/L_0, 8/9.5/1.5/L_1, 6/8/2/L_{k-1}, 1.5/3.5/2/L_{k+1}, 0/1.5/1.5/L_K}{
  \filldraw[fill=blue!10,draw=blue!40] (\l,0) rectangle (\r,\h);
}
\filldraw[fill=red!10,draw=red!50] (3.5,0) rectangle (6,2.5);
\draw[<->,red!75!black] (3.52,1.15) -- (5.98,1.15) node[pos=0.27,above,font=\small] {$P(k)$};
\draw[->] (0,0) -- (11.3,0) node[below right] {$p_S$};
\draw[->] (0,0) -- (0,3.4) node[above] {$\Pop(p_S)$};
\draw[very thick] (0,0) -- (0.75,1.5) -- (2.5,2) -- (4.75,2.5) -- (7,2) -- (8.75,1.5) -- (10,1) -- (10.5,0);
\draw[line width=2.2pt,red!75!black] (4.75,2.5) -- (7,2);
\foreach \x/\y in {0.75/1.5,2.5/2,8.75/1.5,10/1}
  \fill (\x,\y) circle (2.3pt);
\foreach \x/\y in {4.75/2.5,7/2}
  \fill[red!75!black] (\x,\y) circle (2.6pt);
\foreach \x in {0,10.5}
  \filldraw[fill=white,thick] (\x,0) circle (2.3pt);
\draw[densely dashed,red!75!black] (4.75,2.5) -- (4.75,0) (7,2) -- (7,0);
\draw[densely dotted] (0,2.5) -- (4.75,2.5) (0,2) -- (2.5,2);
\draw[densely dotted,gray] (2.5,2) -- (2.5,0) (8.75,1.5) -- (8.75,0) (0.75,1.5) -- (0.75,0) (10,1) -- (10,0);
\foreach \c/\lab in {10/L_0,8.75/L_1,7/L_{k-1},2.5/L_{k+1},0.75/L_K}
  \node[font=\scriptsize,text=blue!60!black,fill=blue!10,inner sep=1.5pt] at (\c,0.22) {$\lab$};
\node[font=\scriptsize,text=red!70!black,fill=red!10,inner sep=1.5pt] at (4.75,0.22) {$L_k$};
\node[left,font=\small] at (0,2.5) {$P(k)$};
\node[left,font=\small] at (0,2) {$P(k-1)$};
\foreach \x/\lab in {0/0,4.75/N-m_k,7/N-m_{k-1},10.5/N}
  \node[below,font=\small] at (\x,-0.03) {$\lab$};
\draw[densely dotted] (0.75,0) -- (0.75,-0.5) node[below,font=\small] {$\tfrac{P(K)}{2}$};
\draw[densely dotted] (10,0) -- (10,-0.5) node[below,font=\small] {$N-\tfrac{P(0)}{2}$};
\node[font=\small] (A) at (0.9,2.95) {slope $2$: $2p_S$};
\draw[->] (A.south) -- (0.4,0.85);
\node[font=\small] (B) at (9.9,2.95) {slope $-2$: $2(N-p_S)$};
\draw[->] (B.south) -- (10.25,0.55);
\draw[->,gray] (8.2,3.5) -- (2.8,3.5) node[midway,above,font=\small] {epidemic progress};
\begin{scope}[font=\scriptsize]
\foreach \l/\r/\lab in {0/0.75/2p_S, 0.75/2.5/j=K, 2.5/4.75/j=k+1, 7/8.75/j=k-1, 8.75/10/j=1}
  \draw[decorate,decoration={brace,mirror,amplitude=3pt}] (\l+0.03,-1.25) -- (\r-0.03,-1.25) node[midway,below=3pt] {$\lab$};
\draw[decorate,decoration={brace,mirror,amplitude=3pt}] (10.03,-1.25) -- (10.47,-1.25) node[midway,below=3pt,anchor=north west,xshift=-6pt] {$2(N-p_S)$};
\draw[decorate,decoration={brace,mirror,amplitude=3pt},red!75!black,thick] (4.78,-1.25) -- (6.97,-1.25) node[midway,below=3pt,text=red!75!black] {$j=k$};
\end{scope}
\end{tikzpicture}
\caption{Interpolated function $\Pop(p_S)$ from Eq.~\eqref{eq:Pop}. Each square corresponds to a layer $L_j$, with side $P(j)$. The nodes $(N-m_j,P(j))$ are joined linearly, with $\Pop(0)=\Pop(N)=0$. Braces indicate the intervals of the definition.}
\label{fig:Pop}
\end{figure}
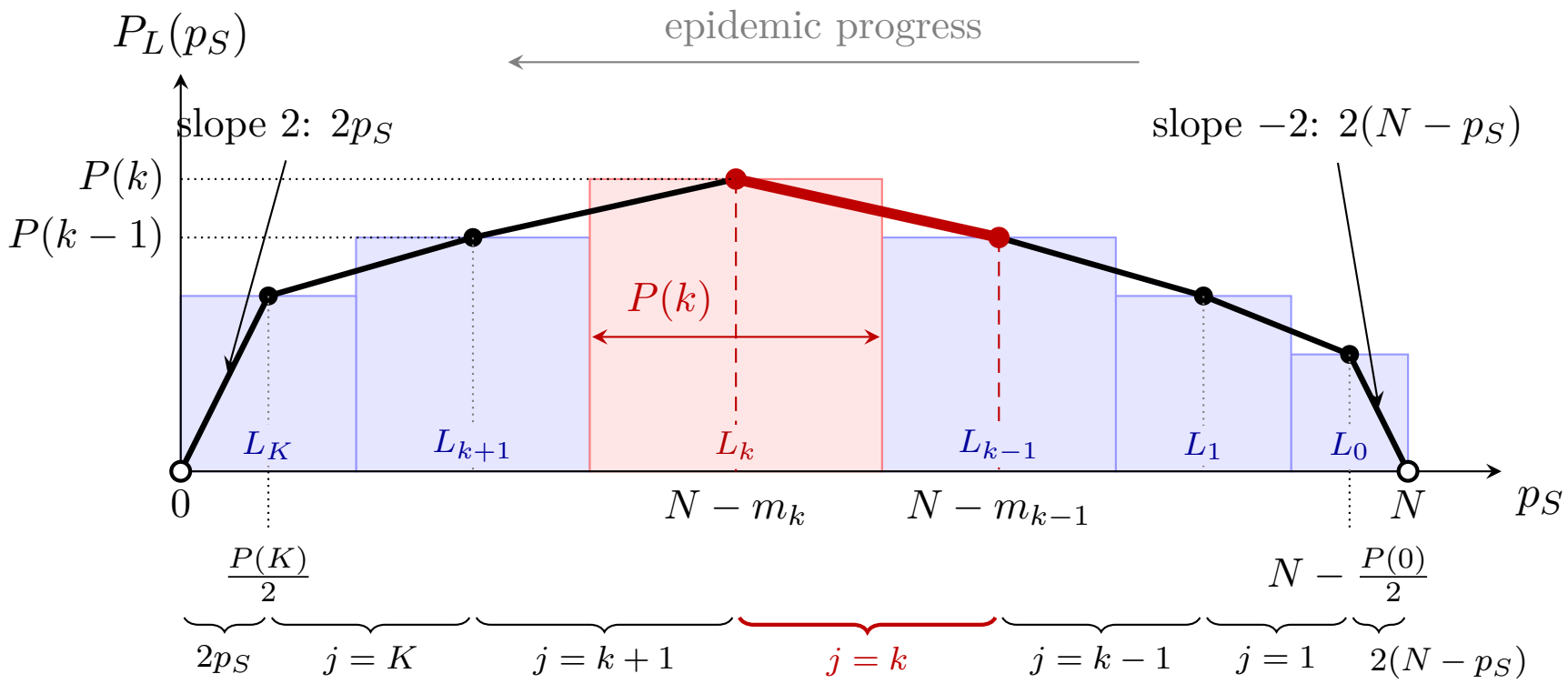

To extend the domain to all $\pS \in [0, N]$, we apply piecewise linear interpolation between the discrete points together with the boundary conditions $\Pop(0) = \Pop(N) = 0$, arriving at
\begin{equation}\label{eq:Pop}
\Pop(p_S) =
\begin{cases}
2(N - p_S), & \text{if } N - \dfrac{P(0)}{2} \le p_S \le N, \\[6pt]
P(k-1) + 2\dfrac{P(k) - P(k-1)}{P(k) + P(k-1)} \Big(N - p_S - m_{k-1}\Big), & \text{if } N - m_k \le p_S \le N - m_{k-1} \quad (1 \le k \le K), \\[8pt]
2 p_S, & \text{if } 0 \le p_S \le \dfrac{P(K)}{2}.
\end{cases}
\end{equation}
Figure~\ref{fig:Pop} illustrates the resulting function $\Pop(p_S)$ defined in Eq.~\eqref{eq:Pop}.

\subsubsection*{Closed Formula for $\Pop(p_S)$}

When $P(k)$ has a closed-form expression, we can derive explicit expressions for $P_c(k)$ and for the midpoints $m_k$. Then, in Eq.~\eqref{eq:Popdisc}, we can obtain a function $\K(\pS)$ such that the condition $\pS = N - m_k$ implies $k = \K(\pS)$, arriving at the closed-form expression
\begin{equation*}
\Pop(\pS) = P(\K(\pS)).
\end{equation*}
Consider, for instance, a regular rectangular grid of $M_1 \times M_2$ districts (with $M_1 \geq M_2$) where each interior district is connected to its four neighbors. When the layer $L_0$ consists of a district located in a corner, $P(k)$ takes the form
\begin{equation} \label{eq:Pk}
P(k) = N_d \min(k+1, M_2, M_1 + M_2 - 1 - k),
\end{equation}
which describes a symmetric trapezoid. In that case, straightforward calculations (see Appendix~\ref{sec:closed_form} for details) lead to the closed formula
\begin{equation} \label{eq:PopRect}
\Pop(\pS) = \min\!\left(\sqrt{2N_d(N-\pS)}, \, N_d M_2, \, \sqrt{2N_d\pS}\right).
\end{equation}

\subsection{Derivation of $\dot{p}_R$}

To derive an expression for $\dot{p}_R(t)$, we assume that every district remains active (i.e., in compartment $\pI$) for a fixed duration $\tau$ after being reached by the epidemic.

\begin{figure}[ht]
\centering
\begin{tikzpicture}[
  >=stealth, scale=0.8, font=\small,
  midarrow/.style={postaction={decorate}, decoration={markings, mark=at position .3 with {\arrow{stealth}}}}
]
\def\T{2}
\def\L{8}
\fill[blue!10] (0,0) -- (\L,0) -- (\L,\T) -- (\T,\T) -- cycle;
\foreach \s in {1.5,3,4.5,6}
  \draw[blue!45, midarrow] (\s,0) -- (\s+\T,\T);
\node[fill=blue!10, inner sep=1.5pt, text=blue!60!black] at (5.25,1) {$\rho(t,a)=-\dot{p}_S(t-a)$};
\draw[->] (0,0) -- (\L+0.6,0) node[right] {$t$};
\draw[->] (0,0) -- (0,\T+1.2) node[above] {$a$};
\draw[thick] (0,0) -- (\L,0);
\draw[thick] (0,0) -- (\T,\T) -- (\L,\T);
\draw[dashed, gray] (\T,\T) -- (\T+0.9,\T+0.9) node[above right, black] {$a=t$};
\draw[dotted, gray] (0,\T) -- (\T,\T);
\draw[dotted, gray] (\T,0) -- (\T,\T);
\draw (0.08,\T) -- (-0.08,\T) node[left] {$\tau$};
\draw (\T,0.08) -- (\T,-0.08) node[below] {$\tau$};
\node at (0.62,1.45) {$\rho=0$};
\node at (5.6,2.6) {$\rho=0$ \quad (age $a>\tau$)};
\end{tikzpicture}
\caption{Domain of $\rho(t,a)$. The function is nonzero only in the shaded region $0\le a\le\tau$, $t\ge a$, and it is constant along the characteristics $t-a=\text{const}$ (thin lines). Each characteristic carries the value $-\dot{p}_S(t-a)$ from the boundary $a=0$ to $a=\tau$.}
\label{fig:rho}
\end{figure}
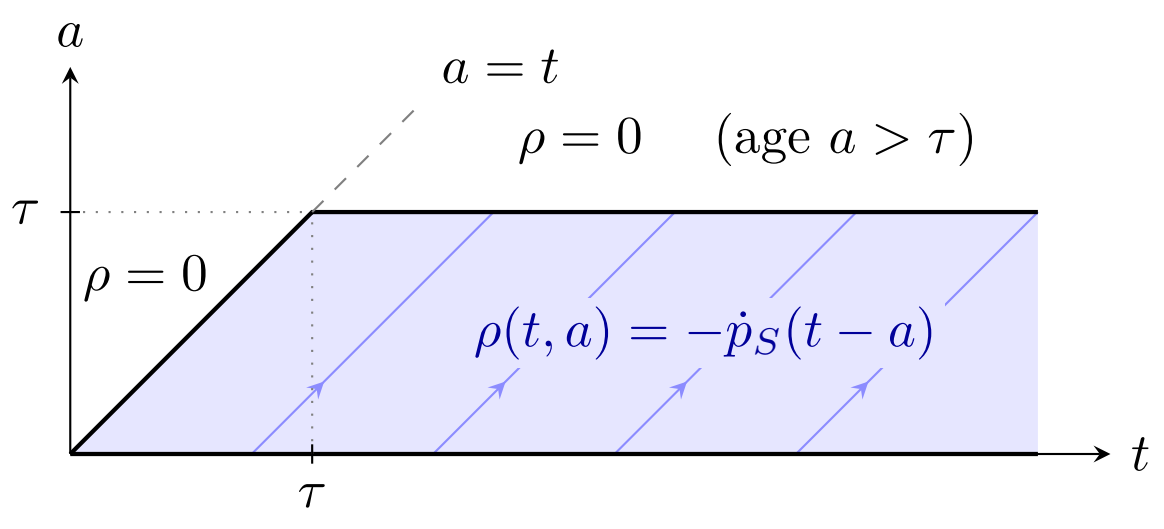

Let $\rho(t,a)$ denote the rate at which individuals enter compartment $p_I$ at time $t - a$, where $a \in [0, \tau]$ is the time elapsed since the entrance of the corresponding district into $p_I$ (that is, $a$ can be viewed as the district infection age). Since individuals leave $\pS$ and enter $\pI$ at rate $-\dpS$, we define
\begin{equation}
\rho(t,a) \coloneq
\begin{cases}
-\dpS(t-a), & \text{if } 0 \le a \le \tau \text{ and } t \ge a, \\[2pt]
0, & \text{otherwise, i.e., if } a > \tau \text{ or } t < a.
\end{cases}
\end{equation}
See Figure~\ref{fig:rho} for an illustration of the domain of $\rho$. The boundary values have a natural interpretation:
\begin{align*}
\rho(t,0) &= -\dpS(t), && \text{for } t \ge 0 \quad \text{(rate of new district entries at time } t\text{)}, \\
\rho(t,\tau) &= -\dpS(t-\tau) = \dpR(t), && \text{for } t \ge \tau \quad \text{(rate of district exits at time } t\text{)}.
\end{align*}
Since districts remain in $\pI$ only up to infection age $a = \tau$ before exiting into $\pR$, the active population in $\pI$ (for $t \ge \tau$, and up to the initial seed $p_I(0)$ for $0 \le t < \tau$) is given by integrating over the active age window $[0, \tau]$:
\begin{equation} \label{eq:pI}
p_I(t) = \int_{0}^{\tau} \rho(t,a)\,\diff a.
\end{equation}

\begin{figure}[htb]
\centering
\begin{tikzpicture}[scale=0.7]
\begin{axis}[
  width=9cm, height=6cm,
  axis lines=left,
  xmin=0, xmax=1.12,
  ymin=0, ymax=0.85,
  xtick={0,1},
  xticklabels={$0$, $\tau$},
  ytick={0.7},
  yticklabels={$R_\infty$},
  xlabel={$a$},
  ylabel={$r(a)$},
  every axis x label/.style={at={(ticklabel* cs:1)},anchor=west},
  every axis y label/.style={at={(ticklabel* cs:1)},anchor=south},
]
\addplot[thick, blue, domain=0:1, samples=100] {0.7*(1-cos(deg(pi*x)))/2};
\addplot[dashed, gray] coordinates {(0,0.7) (1,0.7) (1,0)};
\addplot[only marks, mark=*, mark size=1.5pt] coordinates {(0,0) (1,0.7)};
\end{axis}
\end{tikzpicture}
\caption{A typical profile for the function $r(a)$.}
\label{Fig:r}
\end{figure}
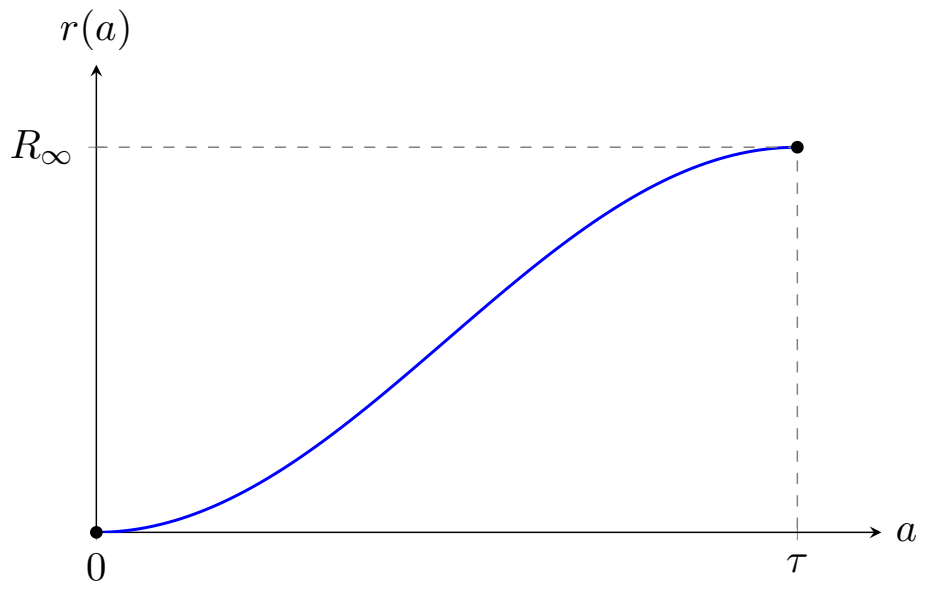

For $a \in [0, \tau]$, let $r(a)$ denote the fraction of the population in a district that has recovered $a$ units of time after the district was reached by the outbreak. Since we assume that all districts follow the same dynamics, $r(a)$ is a fixed profile describing the fraction of recovered individuals as a function of the age of the outbreak. The total recovered population in $p_I$ at time $t$ can then be computed as
\begin{equation} \label{eq:R}
R(t) = \int_{0}^{\tau} r(a)\rho(t,a)\,\diff a.
\end{equation}
Notice also that
\begin{equation}\label{eq:rbounds}
r(0) = 0, \quad r(\tau) = R_{\infty},
\end{equation}
where $R_{\infty}$ represents the fraction of recovered individuals in each district after the epidemic finishes. Figure~\ref{Fig:r} shows a typical profile for $r(a)$.

To deduce a simple expression for $\dot{p}_R(t) = \rho(t, \tau)$, we approximate $r(a)$ by the linear law
\begin{equation} \label{eq:r}
r(a) = R_{\infty} \left(\frac{a}{\tau}\right).
\end{equation}
Using Eq.~\eqref{eq:r} in Eq.~\eqref{eq:R}, the number of recovered individuals in the population $p_I(t)$ becomes
\begin{equation} \label{eq:Rlinear}
R(t) = \int_0^\tau r(a) \rho(t,a)\,\diff a = \frac{R_{\infty}}{\tau} \int_0^\tau a \rho(t,a)\,\diff a.
\end{equation}

We also assume that $\rho(t,a)$ follows a piecewise linear law, considering the two regimes described below. First, we derive the dynamics of $p_R$ during the early stages of the epidemic, when $S$ is large and $R$ is small. We then deduce its dynamics during the later stages of the outbreak, when the epidemic has fully propagated.

\subsubsection{Case 1: Early Outbreak --- Small $R$}

At the beginning of an epidemic wave, infected districts have not yet had sufficient time to reach the maximum infection duration $\tau$. Consequently, the density $\rho(t,a)$ vanishes at some age $a_{\mathrm{max}} \leq \tau$, so that $\rho(t,\tau)$ is identically zero. In this regime, we approximate the density using the piecewise linear expression
\begin{equation} \label{Eq:profilerho}
\rho(t,a) =
\begin{cases}
\rho(t,0)\left(1 - \dfrac{a}{a_{\mathrm{max}}}\right), & \text{for } 0 \leq a \leq a_{\mathrm{max}} \text{ and } t \geq a, \\[6pt]
0, & \text{for } a_{\mathrm{max}} < a \leq \tau.
\end{cases}
\end{equation}
Substituting Eq.~\eqref{Eq:profilerho} into Eq.~\eqref{eq:pI} yields
\begin{equation*}
p_I(t) = \rho(t,0) \int_0^{a_{\mathrm{max}}} \left(1 - \frac{a}{a_{\mathrm{max}}}\right) \diff a = \frac{1}{2} a_{\mathrm{max}} \rho(t,0),
\end{equation*}
which implies
\begin{equation}\label{rho0}
\rho(t,0) = \frac{2 p_I(t)}{a_{\mathrm{max}}}.
\end{equation}
From Eqs.~\eqref{eq:Rlinear} and~\eqref{Eq:profilerho}, we obtain
\begin{equation}\label{Rtamax}
R(t) = \frac{R_{\infty}}{\tau} \int_0^{a_{\mathrm{max}}} a \cdot \rho(t,0) \left(1 - \frac{a}{a_{\mathrm{max}}}\right) \diff a = \frac{R_{\infty} a_{\mathrm{max}}^2 \rho(t,0)}{6 \tau}.
\end{equation}
Using Eqs.~\eqref{rho0} and~\eqref{Rtamax} together with $a_{\max} \le \tau$ and $p_I(t) = R(t) + S(t) + I(t)$, we arrive at
\begin{equation}\label{eq:regime}
R(t) \le \frac{R_\infty}{3 - R_\infty}\bigl(S(t) + I(t)\bigr).
\end{equation}
Hence, $\rho$ satisfies Eq.~\eqref{Eq:profilerho} whenever condition~\eqref{eq:regime} holds.

\subsubsection{Case 2: Late Outbreak}

When the epidemic has fully populated the age domain, $\rho(t,a)$ remains strictly positive across the entire interval $a \in [0, \tau]$. We then approximate it using the affine function
\begin{equation} \label{eq:rholinear}
\rho(t,a) = \rho(t,0)\left(1 - \frac{a}{\tau}\right) + \rho(t,\tau)\left(\frac{a}{\tau}\right).
\end{equation}
Substituting Eq.~\eqref{eq:rholinear} into Eq.~\eqref{eq:pI} gives $p_I(t) = \frac{\tau}{2}\left[\rho(t,0) + \rho(t,\tau)\right]$, which implies
\begin{equation} \label{eq:rho0}
\rho(t,0) = \frac{2}{\tau}p_I(t) - \rho(t,\tau).
\end{equation}
Substituting Eq.~\eqref{eq:rholinear} into Eq.~\eqref{eq:Rlinear} and using Eq.~\eqref{eq:rho0} leads to
\begin{equation*}
R(t) = \frac{R_{\infty}}{3}\left[p_I(t) + \frac{\tau}{2}\rho(t,\tau)\right],
\end{equation*}
from which we isolate $\rho(t,\tau) = \frac{6}{\tau R_{\infty}}R(t) - \frac{2}{\tau}p_I(t)$. Using Eq.~\eqref{eq:pIdef}, we obtain
\begin{equation*}
\rho(t,\tau) = \left(\frac{6 - 2R_{\infty}}{\tau R_{\infty}}\right) R(t) - \frac{2}{\tau}\bigl(S(t) + I(t)\bigr).
\end{equation*}
This can be written in the compact form
\begin{equation*}
\rho(t,\tau) = b_1\bigl(R(t) - b_2 S(t) - b_2 I(t)\bigr),
\end{equation*}
with
\begin{equation}\label{eq:b1b2}
b_1 \coloneq \frac{6 - 2R_{\infty}}{\tau R_{\infty}}, \qquad b_2 \coloneq \frac{R_{\infty}}{3 - R_{\infty}}.
\end{equation}

\subsubsection{The Unified Expression for $\dot{p}_R$}

Combining Cases 1 and 2 above, and recalling that $\dot{p}_R(t) = \rho(t,\tau)$, we obtain
\begin{equation} \label{eq:dotpR}
\dot{p}_R(t) =
\begin{cases}
0, & \text{if } R(t) \le b_2\bigl(S(t) + I(t)\bigr) \text{ or } S(t) = 0, \\[4pt]
b_1\bigl(R(t) - b_2(S(t) + I(t))\bigr), & \text{otherwise.}
\end{cases}
\end{equation}
The condition $S(t) = 0$ in the first branch ensures that $\dot{p}_R(t) = 0$ whenever no susceptible individuals are present, thereby preventing any removal of individuals from $p_I$ under that condition. Observe in the second branch of Eq.~\eqref{eq:dotpR} that the larger $R(t)$ is relative to $S(t) + I(t)$, the faster $p_R$ grows: when $R(t)$ is large, a larger share of districts are completing their outbreak and moving into $p_R$.

\subsection{The Complete pSIR Model}\label{sec:complete}

From Eqs.~\eqref{eq:SIRp}, \eqref{eq:dpS}, and~\eqref{eq:dotpR}, the complete pSIR model is given by
\begin{equation} \label{eq:pSIR}
\begin{split}
\dot{p}_S(t) &= -\frac{\Pop(p_S)}{\Delta t} \min\!\left(1, (I/I_{\mathrm{peak}})^2\right), \\
\dot{S}(t) &= -\dpS(t) - \beta \frac{S(t) I(t)}{S(t)+I(t)+R(t)} - \Sinf \dpR(t), \\
\dot{I}(t) &= \beta \frac{S(t) I(t)}{S(t)+I(t)+R(t)} - \gamma I(t), \\
\dot{R}(t) &= \gamma I(t) - \Rinf \dpR(t), \\
\dot{p}_R(t) &=
\begin{cases}
0, & \text{if } R(t) \le b_2\bigl(S(t)+I(t)\bigr) \text{ or } S(t) = 0, \\[2pt]
b_1\bigl(R(t) - b_2(S(t) + I(t))\bigr), & \text{otherwise.}
\end{cases}
\end{split}
\end{equation}

Besides the classical SIR parameters $\beta$, $\gamma$, and $N$, the model uses three additional parameters: the number of districts $D$, the mean mobility $p$, and the mean connectivity $c$. From these six basic parameters, the following dependent parameters are determined:
\begin{itemize}
\item $\Rinf$ and $\Sinf = 1 - \Rinf$ (the final fractions of recovered and susceptible individuals) satisfy the implicit relation
\begin{equation*}
\Rinf = 1 - e^{-\mathcal{R}_0 \Rinf}, \qquad \text{with } \mathcal{R}_0 = \frac{\beta}{\gamma};
\end{equation*}
\item $\Delta t$ (the time delay between consecutive layers) is computed from Eq.~\eqref{eq:dt};
\item $I_{\mathrm{peak}}$ (the peak value of infected individuals in a district) is obtained from Eq.~\eqref{eq:Ipeak};
\item $b_1$ and $b_2$ are computed from Eq.~\eqref{eq:b1b2}, which also depends on $\tau$ (the duration of the outbreak in a district). This parameter can be estimated as
\begin{equation} \label{eq:tau}
\tau = \frac{2\log(20)}{\beta - \gamma},
\end{equation}
approximating the time elapsed from when the number of infections $I_i(t)$ rises above $5\%$ of the peak $I_{\mathrm{peak}}$ until it falls back below $5\%$ of $I_{\mathrm{peak}}$.
\end{itemize}
In addition, the model involves the function $\Pop(p_S)$, which depends on the network topology and can be computed in general from the layer populations $P(k)$ using Eq.~\eqref{eq:Pop}. Alternatively, closed formulas such as Eq.~\eqref{eq:PopRect} can be derived for regular topologies. Thus, the complete pSIR model with all its parameters is fully characterized from the parameters and topology of the metapopulation model that it approximates.

\section{Properties of the pSIR Model}\label{sec:properties}

In this section, we establish the fundamental mathematical and epidemiological properties of the pSIR system~\eqref{eq:pSIR}, including mass conservation, positivity of the active population, forward invariance, asymptotic convergence, final epidemic size relations, and the calculation of the basic reproduction number $\Rz$.

Throughout this section, we assume that the epidemiological and spatial parameters satisfy $\beta, \gamma, \Delta t, I_{\mathrm{peak}}, b_1, b_2 > 0$, with $\Rinf \in (0, 1)$ and $\Sinf = 1 - \Rinf \in (0, 1)$. Moreover, we assume that the topological activation function $\Pop : [0, N] \to [0, \infty)$ is continuous (and locally Lipschitz on $(0, N]$) with $\Pop(0) = 0$. Recall that the active population is given by $\pI(t) \coloneq S(t) + I(t) + R(t)$, and define the invariant state simplex
\begin{equation}\label{eq:simplex}
\Delta_N \coloneq \left\{ (p_S, S, I, R, p_R) \in \mathbb{R}_+^5 : p_S + S + I + R + p_R = N \right\}.
\end{equation}

\subsection{Well-Posedness and Forward Invariance}

We first verify that total population is conserved, that the active population $\pI(t)$ remains strictly positive in finite time (ensuring that the standard-incidence denominator never vanishes), and that trajectories remain confined to $\Delta_N$.

\begin{lemma}[Conservation of total population]\label{lem:conservation}
Along any solution of system~\eqref{eq:pSIR}, the total population $p_S(t) + S(t) + I(t) + R(t) + p_R(t) = N$ is conserved for all $t \ge 0$, and the active population satisfies
\begin{equation}\label{eq:dpI}
\dot{p}_I(t) = -\dpS(t) - \dpR(t).
\end{equation}
\end{lemma}

\begin{proof}
This follows immediately by summing the five equations in~\eqref{eq:pSIR} and using $\Sinf + \Rinf = 1$.
\end{proof}

\begin{lemma}[Non-vanishing $p_I$]\label{lem:pI_positive}
Let $x(t) = (p_S(t), S(t), I(t), R(t), p_R(t)) \in \Delta_N$ be a solution of~\eqref{eq:pSIR} on $[0, T)$ with $\pI(0) > 0$. Then
\begin{equation}\label{eq:pI_lower_bound}
\pI(t) \ge \pI(0) e^{-b_1 t} > 0 \qquad \text{for all } t \in [0, T).
\end{equation}
In particular, $\pI(t)$ never vanishes in finite time.
\end{lemma}

\begin{proof}
Since $S(t), I(t), R(t) \ge 0$ and $b_1, b_2 > 0$, we have $R(t) - b_2(S(t) + I(t)) \le R(t) \le \pI(t)$, so both branches of the $\dpR(t)$ equation in~\eqref{eq:pSIR} satisfy
\begin{equation}\label{eq:dpR_bound}
0 \le \dpR(t) \le b_1 \pI(t).
\end{equation}
Furthermore, since $\Pop(p_S) \ge 0$, we have $-\dpS(t) \ge 0$. Substituting these bounds into~\eqref{eq:dpI} yields
\begin{equation*}
\dot{p}_I(t) = -\dpS(t) - \dpR(t) \ge -b_1 \pI(t).
\end{equation*}
Applying Gr\"onwall's inequality gives~\eqref{eq:pI_lower_bound}.
\end{proof}

\begin{proposition}[Well-posedness and forward invariance]\label{prop:invariance}
For every initial condition $x(0) \in \Delta_N$ with $\pI(0) > 0$, system~\eqref{eq:pSIR} admits a global forward solution that remains in $\Delta_N \cap \{\pI > 0\}$ for all $t \ge 0$.
\end{proposition}

\begin{proof}
By Lemma~\ref{lem:pI_positive}, $\pI(t) > 0$ on any existence interval $[0, T)$, so the standard-incidence term $\beta S I / \pI$ remains well-defined and locally Lipschitz on $\Delta_N \cap \{\pI > 0\}$. We check the boundary behavior on each face of $\mathbb{R}_+^5$:
\begin{itemize}
\item If $p_S = 0$, then $\Pop(0) = 0$ implies $\dpS = 0$.
\item If $S = 0$, then $\dpR = 0$ and $\dot{S} = -\dpS \ge 0$.
\item If $I = 0$, then $\dot{I} = 0$.
\item If $R = 0$, then $R \le b_2(S + I)$ implies $\dpR = 0$, so $\dot{R} = \gamma I \ge 0$.
\item If $p_R = 0$, then $\dpR \ge 0$.
\end{itemize}
Hence $\mathbb{R}_+^5$ is forward-invariant. Combined with mass conservation (Lemma~\ref{lem:conservation}), every trajectory remains bounded in the compact set $\Delta_N$, which guarantees global existence for all $t \ge 0$.
\end{proof}

\subsection{Transient Dynamics}\label{subsec:transient}

Next, we characterize the instantaneous condition governing epidemic growth or decay and explain how spatial recruitment ($p_S \to \pI$) and post-outbreak removal ($\pI \to p_R$) generate incidence plateaus and multiple waves.

\begin{proposition}[Threshold condition for outbreak growth]\label{prop:threshold}
Let $\mathcal{R}_0 \coloneq \beta / \gamma$. Along any solution in $\Delta_N \cap \{\pI > 0\}$:
\begin{enumerate}
\item If $\mathcal{R}_0 \le 1$, then $\dot{I}(t) < 0$ whenever $I(t) > 0$, so no epidemic outbreak can occur.
\item If $\mathcal{R}_0 > 1$, then $\dot{I}(t) > 0$ at time $t$ if and only if $I(t) > 0$ and
\begin{equation}\label{eq:growth_cond}
\frac{S(t)}{\pI(t)} > \frac{1}{\mathcal{R}_0}.
\end{equation}
\end{enumerate}
\end{proposition}

\begin{proof}
Factoring the $\dot{I}(t)$ equation in~\eqref{eq:pSIR} gives
\begin{equation*}
\dot{I}(t) = \gamma \left(\mathcal{R}_0 \frac{S(t)}{\pI(t)} - 1\right) I(t).
\end{equation*}
Since $0 \le S(t) \le \pI(t)$ and $\pI(t) > 0$ for all $t \ge 0$, we have $S(t)/\pI(t) \le 1$. If $\mathcal{R}_0 \le 1$, then $\mathcal{R}_0 S(t)/\pI(t) \le 1$ (with strict inequality whenever $I(t) > 0$, as $\pI(t) \ge S(t) + I(t) > S(t)$), yielding $\dot{I}(t) < 0$. If $\mathcal{R}_0 > 1$ and $I(t) > 0$, the sign of $\dot{I}(t)$ is strictly positive if and only if~\eqref{eq:growth_cond} holds.
\end{proof}

\begin{remark}\label{rem:multiwave}
In the classical single-population SIR model, the susceptible fraction $S(t)/N$ is strictly non-increasing, which restricts the incidence curve $I(t)$ to at most a single peak. In contrast, differentiating the active susceptible fraction $S(t)/\pI(t)$ in~\eqref{eq:pSIR} yields
\begin{equation*}
\frac{\diff}{\diff t}\left(\frac{S}{\pI}\right) = -\beta \frac{S I}{\pI^2} - \frac{\dpS}{\pI}\left(1 - \frac{S}{\pI}\right) + \frac{\dpR}{\pI}\left(\frac{S}{\pI} - \Sinf\right).
\end{equation*}
Because $-\dpS(t) \ge 0$ replenishes $S/\pI$ as new districts are activated and $\dpR(t) \ge 0$ removes post-outbreak districts (where $S/\pI$ has dropped toward $\Sinf$), the ratio $S(t)/\pI(t)$ need not decrease monotonically: it can hover near $1/\mathcal{R}_0$ to sustain a prolonged plateau or rise back above $1/\mathcal{R}_0$ after a dip, producing multiple infectious waves.
\end{remark}

\subsection{Asymptotic Behavior}\label{subsec:asymptotic}

Below we analyze the long-time limits of all five state variables as $t \to \infty$ and derive the exact final epidemic size relations.

\begin{proposition}[Asymptotic convergence and disease extinction]\label{prop:convergence}
For every initial condition $x(0) \in \Delta_N$ with $\pI(0) > 0$, the compartment $p_S(t)$ is non-increasing, $p_R(t)$ is non-decreasing, the infected population vanishes asymptotically,
\begin{equation}\label{eq:I_extinction}
\lim_{t \to \infty} I(t) = 0,
\end{equation}
and the remaining state variables converge to finite limits
\begin{equation*}
p_S^\infty \coloneq \lim_{t \to \infty} p_S(t), \quad S^\infty \coloneq \lim_{t \to \infty} S(t), \quad R^\infty \coloneq \lim_{t \to \infty} R(t), \quad p_R^\infty \coloneq \lim_{t \to \infty} p_R(t)
\end{equation*}
(with $\pI^\infty \coloneq \lim_{t \to \infty} \pI(t) = N - p_S^\infty - p_R^\infty$) such that $x^\infty \coloneq (p_S^\infty, S^\infty, 0, R^\infty, p_R^\infty) \in \Delta_N$ is an equilibrium of~\eqref{eq:pSIR}.
\end{proposition}

\begin{proof}
By Proposition~\ref{prop:invariance}, solutions remain in $\Delta_N \cap \{\pI > 0\}$ for all $t \ge 0$. Since $\Pop(p_S) \ge 0$ and $\Delta t > 0$, we have $\dpS(t) \le 0$, whereas $\dpR(t) \ge 0$ by definition. Because $p_S(t), p_R(t) \in [0, N]$ are bounded and monotone, their limits $p_S^\infty$ and $p_R^\infty$ exist, and $\pI^\infty$ exists by Lemma~\ref{lem:conservation}.

To establish convergence of the remaining variables, consider the two auxiliary functions
\begin{equation*}
W(t) \coloneq R(t) + \Rinf p_R(t), \qquad U(t) \coloneq S(t) + p_S(t) + \Sinf p_R(t).
\end{equation*}
Differentiating $W(t)$ and $U(t)$ along solutions of~\eqref{eq:pSIR} yields
\begin{equation}\label{eq:WU_deriv}
\dot{W}(t) = \gamma I(t) \ge 0, \qquad \dot{U}(t) = -\beta \frac{S(t) I(t)}{\pI(t)} \le 0.
\end{equation}
Since $W(t)$ is non-decreasing and bounded above by $(1 + \Rinf)N$, the limit $W^\infty \coloneq \lim_{t \to \infty} W(t)$ exists and
\begin{equation*}
\int_0^\infty I(t)\,\diff t = \frac{W^\infty - W(0)}{\gamma} < \infty.
\end{equation*}
Moreover, since $0 \le S(t) \le \pI(t)$ and $0 \le I(t) \le N$ for all $t \ge 0$, the derivative of $I(t)$ is uniformly bounded:
\begin{equation*}
|\dot{I}(t)| = \left|\beta \frac{S(t)}{\pI(t)} - \gamma\right| I(t) \le \max(\gamma, |\beta - \gamma|) N \qquad \text{for all } t \ge 0.
\end{equation*}
Hence $I(t)$ is globally Lipschitz (and thus uniformly continuous) and integrable on $[0, \infty)$. Consequently, $\lim_{t \to \infty} I(t) = 0$.

Similarly, $U(t) \ge 0$ is non-increasing by~\eqref{eq:WU_deriv}, so $U^\infty \coloneq \lim_{t \to \infty} U(t)$ exists. It follows that
\begin{equation*}
S^\infty \coloneq \lim_{t \to \infty} S(t) = U^\infty - p_S^\infty - \Sinf p_R^\infty, \qquad R^\infty \coloneq \lim_{t \to \infty} R(t) = N - p_S^\infty - S^\infty - p_R^\infty
\end{equation*}
both exist, with $x^\infty \in \Delta_N$. Finally, if $S^\infty = 0$, then $x^\infty$ is an equilibrium on $\{S = 0, I = 0\}$; if $S^\infty > 0$, then $S(t) > 0$ for all sufficiently large $t$, so $\lim_{t \to \infty} \dpR(t) = b_1 \max(0, R^\infty - b_2 S^\infty)$. Since $\int_0^\infty \dpR(t)\,\diff t = p_R^\infty - p_R(0) < \infty$, this limit must be zero, i.e., $R^\infty \le b_2 S^\infty$. In either case, $x^\infty$ is an equilibrium of~\eqref{eq:pSIR}.
\end{proof}

\begin{proposition}[Final epidemic size relation]\label{prop:final_size}
For any initial condition $x(0) \in \Delta_N$ with $\pI(0) > 0$, the cumulative recovery integral and cumulative incidence satisfy
\begin{equation}\label{eq:cum_recovery}
\gamma \int_0^\infty I(t)\,\diff t = \bigl(R^\infty - R(0)\bigr) + \Rinf \bigl(p_R^\infty - p_R(0)\bigr),
\end{equation}
and
\begin{equation}\label{eq:cum_incidence}
\int_0^\infty \beta \frac{S(t) I(t)}{\pI(t)}\,\diff t = \bigl(R^\infty - R(0)\bigr) + \Rinf \bigl(p_R^\infty - p_R(0)\bigr) - I(0).
\end{equation}
In particular, if $R(0) = p_R(0) = 0$ and all reached districts eventually complete their outbreaks ($\pI^\infty = 0$, so $S^\infty = R^\infty = 0$ and $p_R^\infty = N - p_S^\infty$), the total fraction of ever-infected individuals factorizes as
\begin{equation}\label{eq:final_size_factorized}
\frac{\gamma}{N}\int_0^\infty I(t)\,\diff t = \Rinf \left(1 - \frac{p_S^\infty}{N}\right),
\end{equation}
while the cumulative new-infection fraction is $\Rinf\bigl(1 - p_S^\infty/N\bigr) - I(0)/N$.
\end{proposition}

\begin{proof}
From~\eqref{eq:pSIR}, we have $\dot{R}(t) + \Rinf \dpR(t) = \gamma I(t)$ for all $t \ge 0$. Integrating both sides over $[0, \infty)$ and applying Proposition~\ref{prop:convergence} directly yields~\eqref{eq:cum_recovery}. Integrating the $\dot{I}(t)$ equation on $[0, \infty)$ and using~\eqref{eq:I_extinction} gives
\begin{equation*}
\int_0^\infty \beta \frac{S(t) I(t)}{\pI(t)}\,\diff t = \gamma \int_0^\infty I(t)\,\diff t + \lim_{t \to \infty} I(t) - I(0) = \gamma \int_0^\infty I(t)\,\diff t - I(0),
\end{equation*}
which combined with~\eqref{eq:cum_recovery} yields~\eqref{eq:cum_incidence}. Setting $R(0) = p_R(0) = 0$ and $R^\infty = 0$, $p_R^\infty = N - p_S^\infty$ gives~\eqref{eq:final_size_factorized}.
\end{proof}

\subsection{Basic Reproduction Number}

The following theorem establishes that the basic reproduction number $\Rz$ of the pSIR model coincides with that of the metapopulation model it approximates and with those of the individual SIR models of the corresponding network.

\begin{theorem} \label{ThmR0}
The basic reproduction number for model~\eqref{eq:pSIR} is given by
\begin{equation}
\mathcal{R}_0 = \frac{\beta}{\gamma},
\end{equation}
and coincides both with the basic reproduction number of the standard SIR model and with that of the metapopulation model of Eq.~\eqref{eq:SIRi}.
\end{theorem}

The proof of this theorem follows from the computation of $\Rz$ for the pSIR model (see Appendix~\ref{AppR0}) and for the metapopulation model (see Appendix~\ref{AppR0network}).

\section{Numerical Examples}\label{SecExamples}

In this section, we present three numerical examples. In the first one, the aggregated pSIR model~\eqref{eq:pSIR} is used to reproduce the simulation results of a synthetic metapopulation model in which the districts are arranged on a regular grid. In the second example, the synthetic metapopulation model is modified so that the districts have heterogeneous populations and are connected through a randomly generated network. Finally, in the third example, the pSIR model is used to fit real data corresponding to the daily number of COVID-19 cases reported in Argentina during the period from March to November 2020.

The reported results can be reproduced using the Python routines available at \url{https://github.com/CIFASIS/pSIR/}. Those routines can also be used to obtain and simulate pSIR approximations of the aforementioned metapopulation systems under arbitrary sets of parameters.

\subsection{Synthetic Metapopulation Model --- Rectangular Grid}

\begin{figure}[htb]
\centering
\begin{tikzpicture}[
  scale=0.9,
  district/.style={circle, draw, fill=white, minimum size=5.5mm, inner sep=0pt, font=\scriptsize}
]
\foreach \i in {1,...,4}{
  \foreach \j in {1,...,5}{
    \coordinate (p\i\j) at (\j, -\i);
  }
}
\foreach \i in {1,...,4}{
  \foreach \j in {1,...,4}{
    \pgfmathtruncatemacro{\jn}{\j+1}
    \draw[gray!70] (p\i\j) -- (p\i\jn);
  }
}
\foreach \i in {1,...,3}{
  \foreach \j in {1,...,5}{
    \pgfmathtruncatemacro{\in}{\i+1}
    \draw[gray!70] (p\i\j) -- (p\in\j);
  }
}
\foreach \i in {1,...,4}{
  \foreach \j in {1,...,5}{
    \node[district] at (p\i\j) {};
  }
}
\node[district, fill=blue!20,   draw=blue!70]   at (p23) {$i,j$};
\node[district, fill=orange!25, draw=orange!85] at (p25) {};
\node[district, fill=red!20,    draw=red!70]    at (p11) {};
\foreach \i/\lab in {1/1, 2/2, 3/3, 4/M_1}
  \node[font=\scriptsize, left=2mm] at (p\i1) {$\lab$};
\node[font=\scriptsize] at (0.45,-3.45) {$\vdots$};
\foreach \j/\lab in {1/1, 2/2, 3/3, 4/4, 5/M_2}
  \node[font=\scriptsize, above=2mm] at (p1\j) {$\lab$};
\node[font=\scriptsize] at (4.5,-0.45) {$\cdots$};
\node[font=\small, left=8mm]  at (1,-2.5) {$i$};
\node[font=\small, above=7mm] at (p13) {$j$};
\begin{scope}[shift={(6.4,-1.8)}, font=\small]
  \node[district, fill=blue!20, draw=blue!70, minimum size=4mm] (a) at (0,0) {};
  \node[right=2mm of a, anchor=west] {interior: 4 neighbors};
  \node[district, fill=orange!25, draw=orange!85, minimum size=4mm] (b) at (0,-0.7) {};
  \node[right=2mm of b, anchor=west] {boundary: 3 neighbors};
  \node[district, fill=red!20, draw=red!70, minimum size=4mm] (c) at (0,-1.4) {};
  \node[right=2mm of c, anchor=west] {corner: 2 neighbors};
\end{scope}
\end{tikzpicture}
\caption{Rectangular grid of $M_1 \times M_2$ districts ($M_1$ rows, $M_2$ columns) with nearest-neighbor connections. Highlighted: an interior district $(i,j)$ and its four neighbors, a boundary district (three), and a corner district (two).}
\label{fig:grid}
\end{figure}
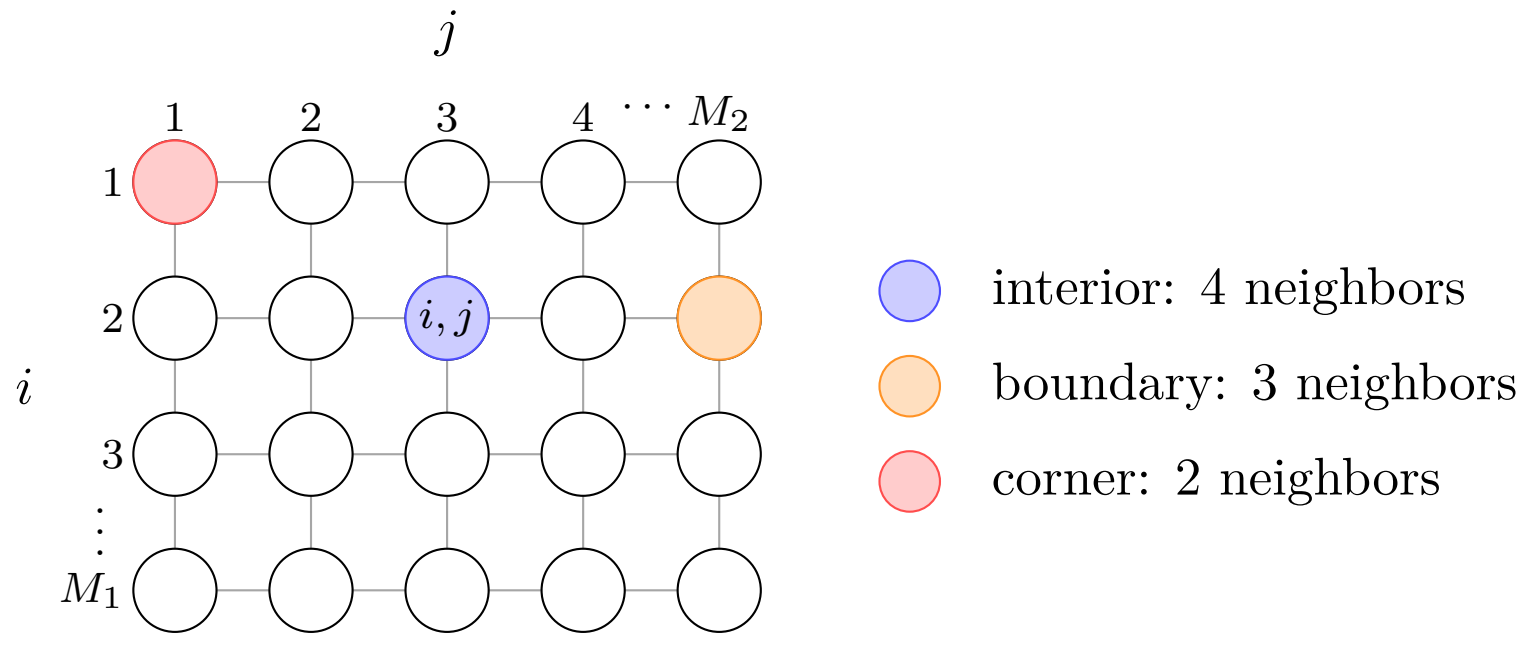

We consider a network consisting of a rectangular grid of $M_1 \times M_2$ districts, indexed by $(i,j)$ with $1 \le i \le M_1$ and $1 \le j \le M_2$. Each district is connected to its nearest neighbors $(i \pm 1, j)$ and $(i, j \pm 1)$ that lie within the grid. Interior districts therefore have four neighbors, boundary districts three, and corner districts two (see Figure~\ref{fig:grid}).

We assume that the $D = M_1 \cdot M_2$ districts have equal population size $N_i = N/D$ for $i = 1, \ldots, D$, and that the dynamics of the districts follow Eq.~\eqref{eq:SIRi}.

For the numerical simulations, we initially consider the parameters $\beta = 1$, $\gamma = 0.5$, $p = 0.005$, and a total population of $N = 10^6$ individuals equally distributed across $D = 20$ districts ($M_1 = 10$, $M_2 = 2$). The connectivity parameters $c_i$ are $2$ for the corner districts, $3$ for the boundary districts, and $4$ for the interior districts. The initial condition corresponds to an epidemic originating with $10$ infected individuals ($I_1(0) = 10$) in a corner district, while the remaining individuals in that district and all individuals living in the other districts are initially susceptible.

The parameters of the pSIR model were obtained following Section~\ref{sec:complete}, resulting in $\Rinf = 0.797$, $\Delta t = 5.229$, $I_{\mathrm{peak}} = 7671.3$, $b_1 = 0.4615$, $b_2 = 0.3617$, and $\tau = 11.983$. The function $\Pop(p_S)$ was computed using the closed form of Eq.~\eqref{eq:PopRect}. The initial conditions for the pSIR model were $I(0) = I_1(0) = 10$, $p_I(0) = N/D = 5 \cdot 10^4$, $S(0) = p_I(0) - I(0)$, $p_S(0) = N - p_I(0)$, and $p_R(0) = R(0) = 0$.

\begin{figure}[htb]
\centering
\includegraphics[width=0.8\linewidth]{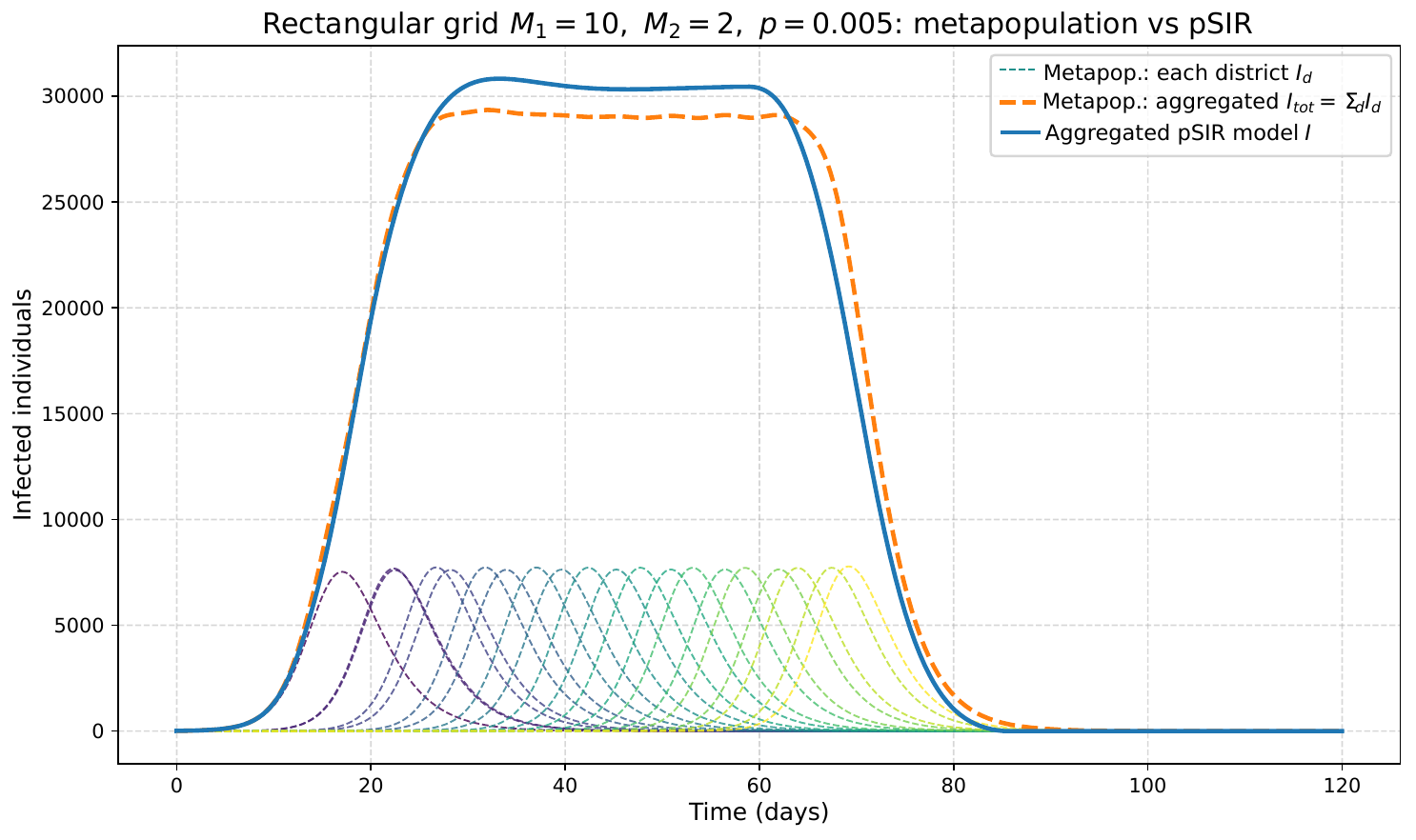}
\caption{Metapopulation vs.\ pSIR aggregated model on a $10 \times 2$ rectangular grid.}
\label{fig:basecase}
\end{figure}

Figure~\ref{fig:basecase} shows the evolution of the infected population in each district, together with the total number of infected individuals across the metapopulation, and compares these results with those obtained using the pSIR aggregated model.

The results show that the aggregated model successfully approximates the characteristic shape and timing of the epidemic dynamics arising from this spatial configuration. In particular, it reproduces the prolonged epidemic plateau followed by a sharp decline, a behavior that standard SIR models cannot capture.

To verify that the pSIR model accurately approximates the metapopulation model under different parameter settings, we repeated the experiment varying the mobility ($p = 0.001$ and $p = 0.05$) for two network configurations: an elongated grid with dimensions $M_1 \times M_2 = 10 \times 2$ and a compact one with $M_1 \times M_2 = 6 \times 4$.

\begin{figure}[htb]
\centering
\includegraphics[width=\linewidth]{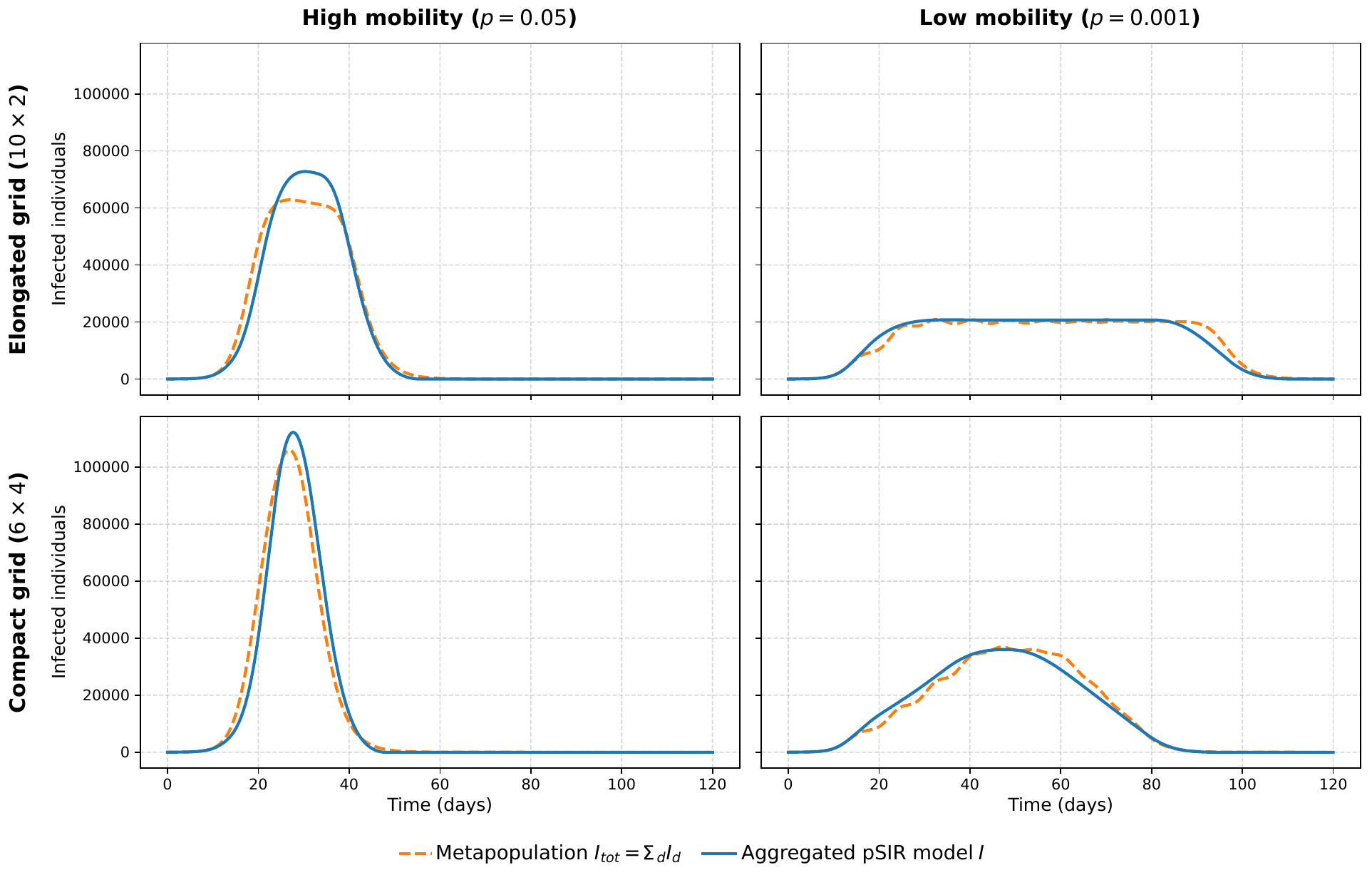}
\caption{Rectangular grid: Scenarios with varying mobility and grid aspect ratio.}
\label{fig:variations}
\end{figure}

Figure~\ref{fig:variations} shows that the pSIR aggregated model successfully captures the epidemic dynamics in all four scenarios. As expected, mobility has a strong influence on the speed of epidemic propagation: higher mobility leads to faster spread across districts and a more pronounced epidemic peak, whereas lower mobility results in slower propagation and a more prolonged plateau.

The spatial configuration also affects the epidemic dynamics. Comparing the elongated and compact grids under the same mobility conditions, we observe that the epidemic spreads more rapidly in the compact configuration. This can be attributed to the geometry of the grid, which allows the infection to reach a larger number of districts in fewer propagation steps. These results demonstrate that the aggregated model captures the effects of both mobility and spatial structure while retaining a substantially simpler representation of the epidemic.

\subsection{Random Spatial Network}

We now consider a metapopulation model composed of $D = 40$ districts whose spatial distribution and mobility relationships are randomly generated. More precisely, the district locations are randomly sampled within a rectangular domain whose shape is controlled by its dimensions $L_x$ and $L_y$. The total population $N$ is then distributed among the $D$ districts, allowing for population heterogeneity.

The network is constructed by establishing mobility links between districts so as to ensure connectivity while imposing a prescribed average number of connections per district $c_{\mathrm{con}}$. Mobility flows are then randomly assigned to connected pairs of districts, taking into account a mean mobility value $p$, their populations, and their spatial distances.

\begin{figure}[htb]
\centering
\includegraphics[width=\linewidth]{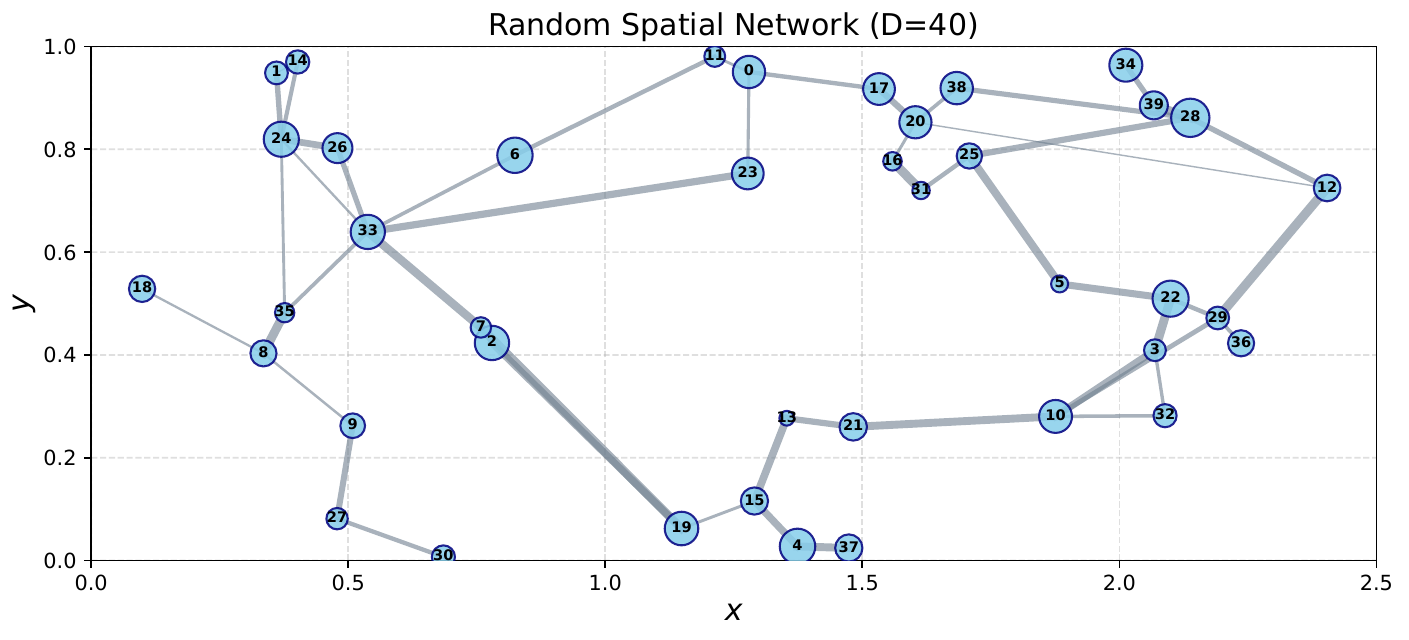}
\caption{Spatial configuration and mobility network: First Scenario.}
\label{fig:rand_pop}
\end{figure}

As a first scenario, we consider a population of $N = 10^6$ individuals heterogeneously distributed among $D = 40$ districts, with parameter values $L_x = 2.5$, $L_y = 1$, $c_{\mathrm{con}} = 2.5$, $p = 0.005$, $\beta = 1$, and $\gamma = 0.5$. As in the previous example, the epidemic is initialized with $I_1(0) = 10$ infected individuals in one of the districts. Figure~\ref{fig:rand_pop} shows the resulting spatial configuration and mobility network.

The network of Figure~\ref{fig:rand_pop} has eight layers with populations $P(0) = 33{,}427$, $P(1) = 78{,}569$, $P(2) = 110{,}849$, $P(3) = 160{,}907$, $P(4) = 203{,}021$, $P(5) = 209{,}545$, $P(6) = 162{,}810$, and $P(7) = 40{,}868$. From these values, the function $\Pop(p_S)$ is given by the interpolation formula of Eq.~\eqref{eq:Pop}.

Proceeding as in the previous example, the pSIR model parameters evaluate to $\Rinf = 0.797$, $\Delta t = 5.394$, $b_1 = 0.4615$, $b_2 = 0.3617$, $\tau = 11.983$, and $I_{\mathrm{peak}} = 5128.7$, where this last value was computed using Eq.~\eqref{eq:Ipeak} with $N_d$ replaced by $P(0)$.

\begin{figure}[htb]
\centering
\includegraphics[width=0.8\linewidth]{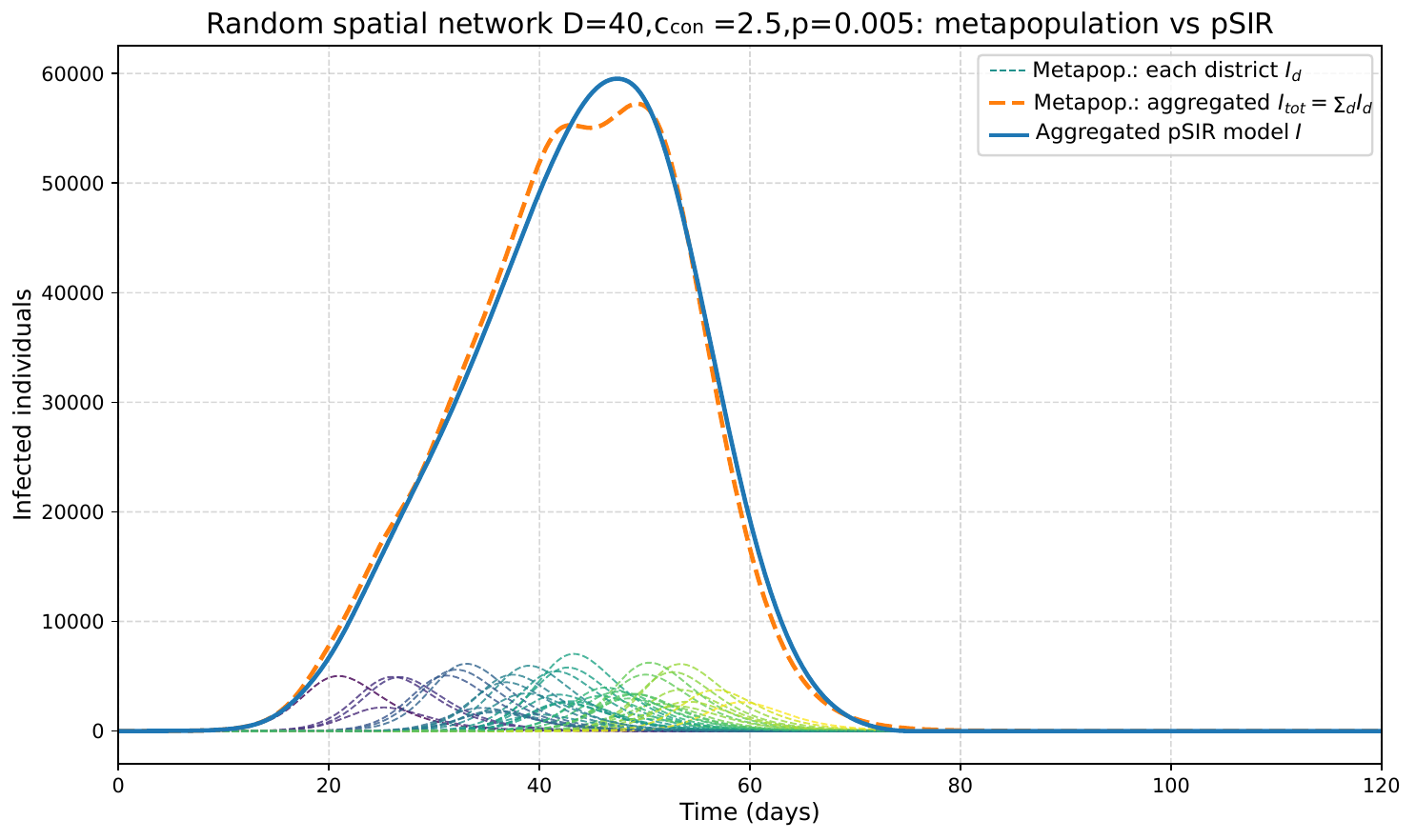}
\caption{Evolution of infected individuals: First Scenario.}
\label{fig:random_base_case}
\end{figure}

Figure~\ref{fig:random_base_case} presents the time evolution of the number of infected individuals in each district, together with the total number of infected individuals in the metapopulation model, compared against the trajectory obtained by the aggregated pSIR model.

\begin{figure}[htb]
\centering
\includegraphics[width=\linewidth]{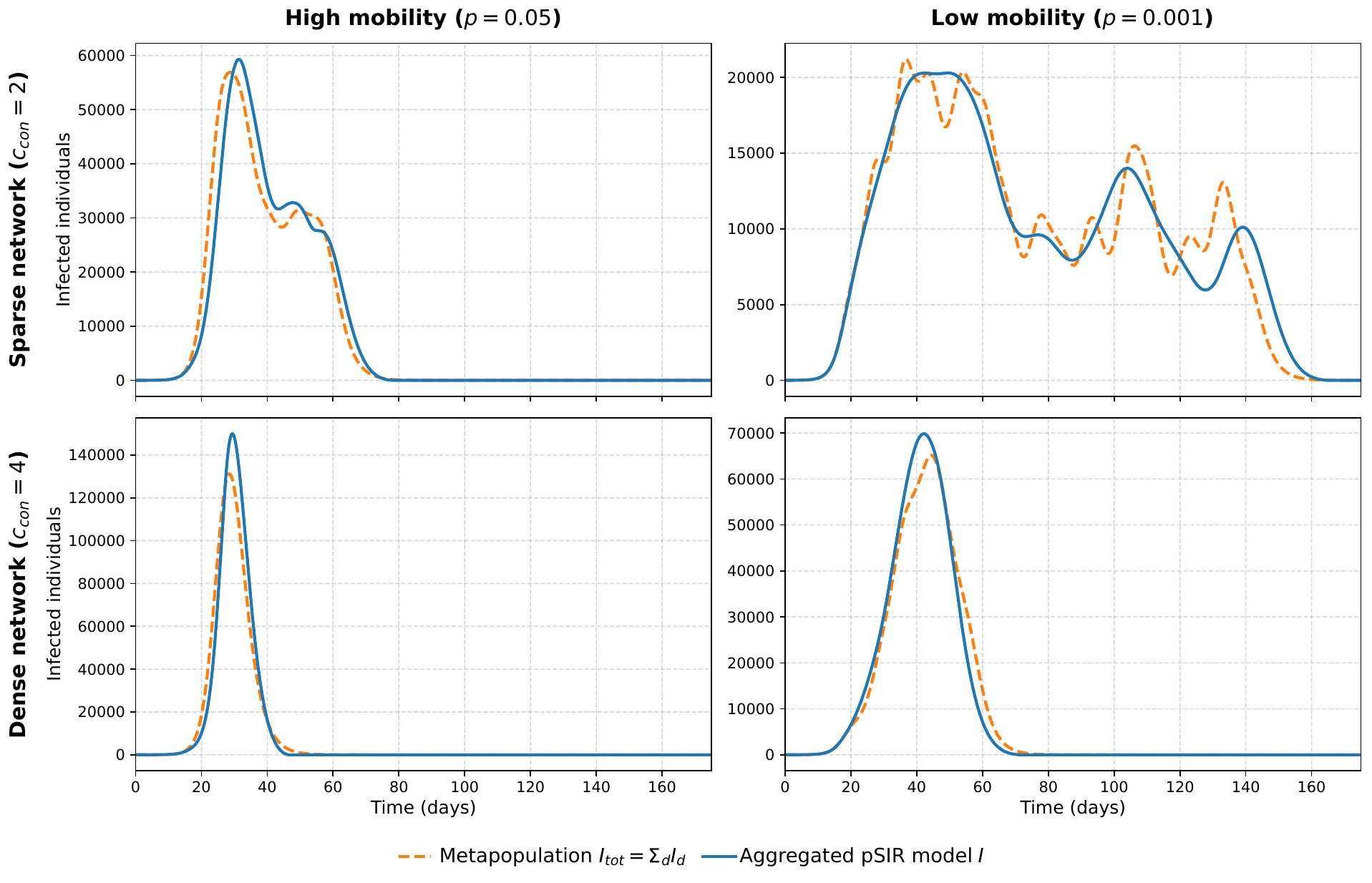}
\caption{Random spatial network: Scenarios with varying mobility and connectivity parameters.}
\label{fig:random_variations}
\end{figure}

We also considered additional scenarios by varying the average mobility level ($p = 0.001$ and $p = 0.05$) and the average network connectivity ($c_{\mathrm{con}} = 2$ and $c_{\mathrm{con}} = 4$). Figure~\ref{fig:random_variations} compares the aggregated number of infected individuals obtained from the metapopulation simulations with the trajectories of the aggregated pSIR model across the four combinations of mobility and connectivity levels.

Once again, the aggregated pSIR model reproduces the metapopulation dynamics remarkably well across all scenarios, capturing the timing of the epidemic, the occurrence of multiple waves, and the overall shape of the infection curves.

\subsection{Fitting Real Data}

\begin{figure}[htbp]
\centering
\includegraphics[width=\linewidth]{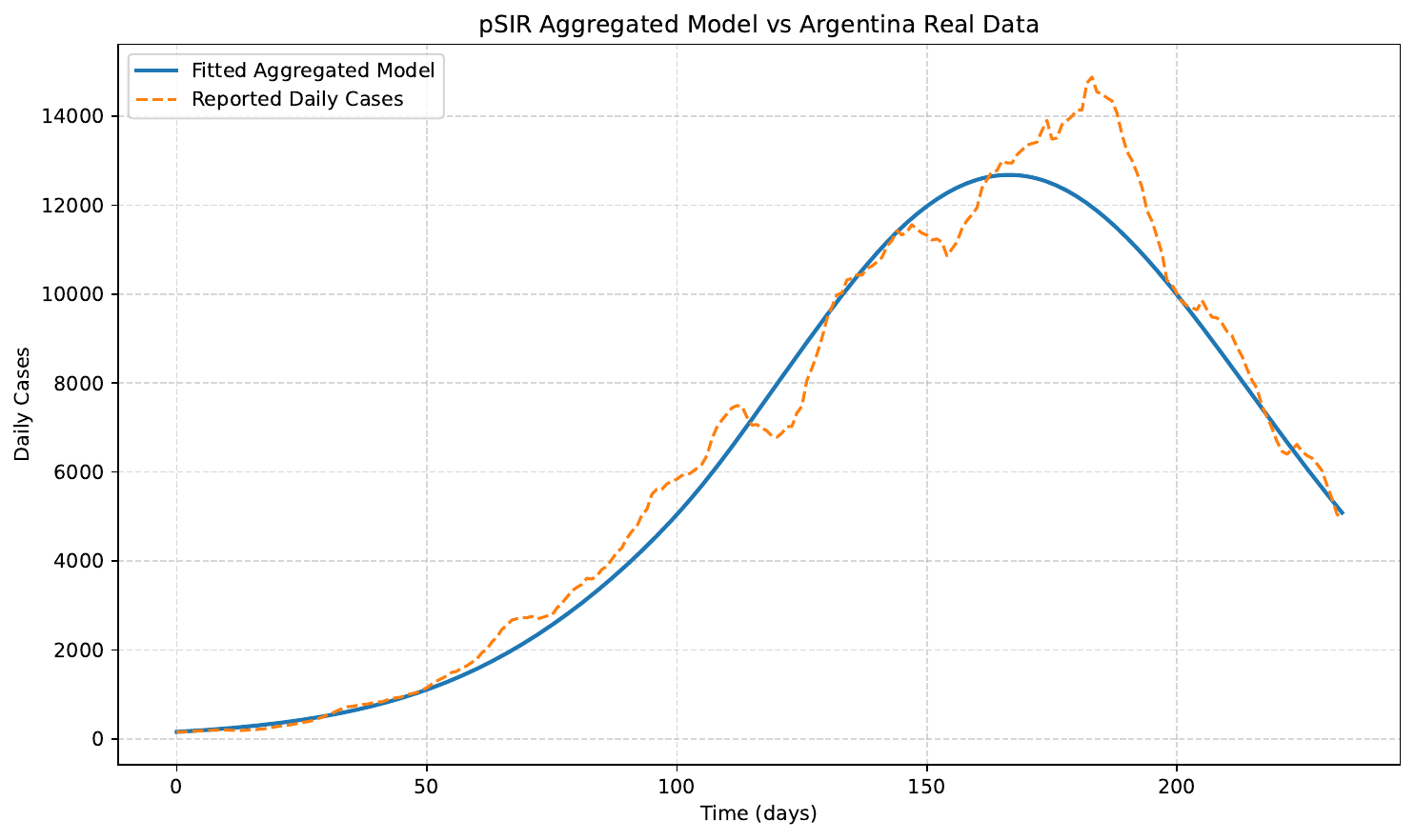}
\caption{Daily cases of COVID-19 in Argentina (with 7-day moving average) and simulation results with the pSIR aggregated model.}
\label{fig:arg_fitted}
\end{figure}

The third example corresponds to the first wave of COVID-19 in Argentina that occurred during the period March--November 2020. The orange curve in Figure~\ref{fig:arg_fitted} shows the daily number of reported cases (processed with a moving average filter of length $W = 7$). The raw data can be downloaded from \url{https://covidstats.com.ar}, and the processed daily cases dataset is available at \url{https://fceia.unr.edu.ar/~kofman/files/nacion.csv}.

\begin{figure}[htbp]
\centering
\includegraphics[width=\linewidth]{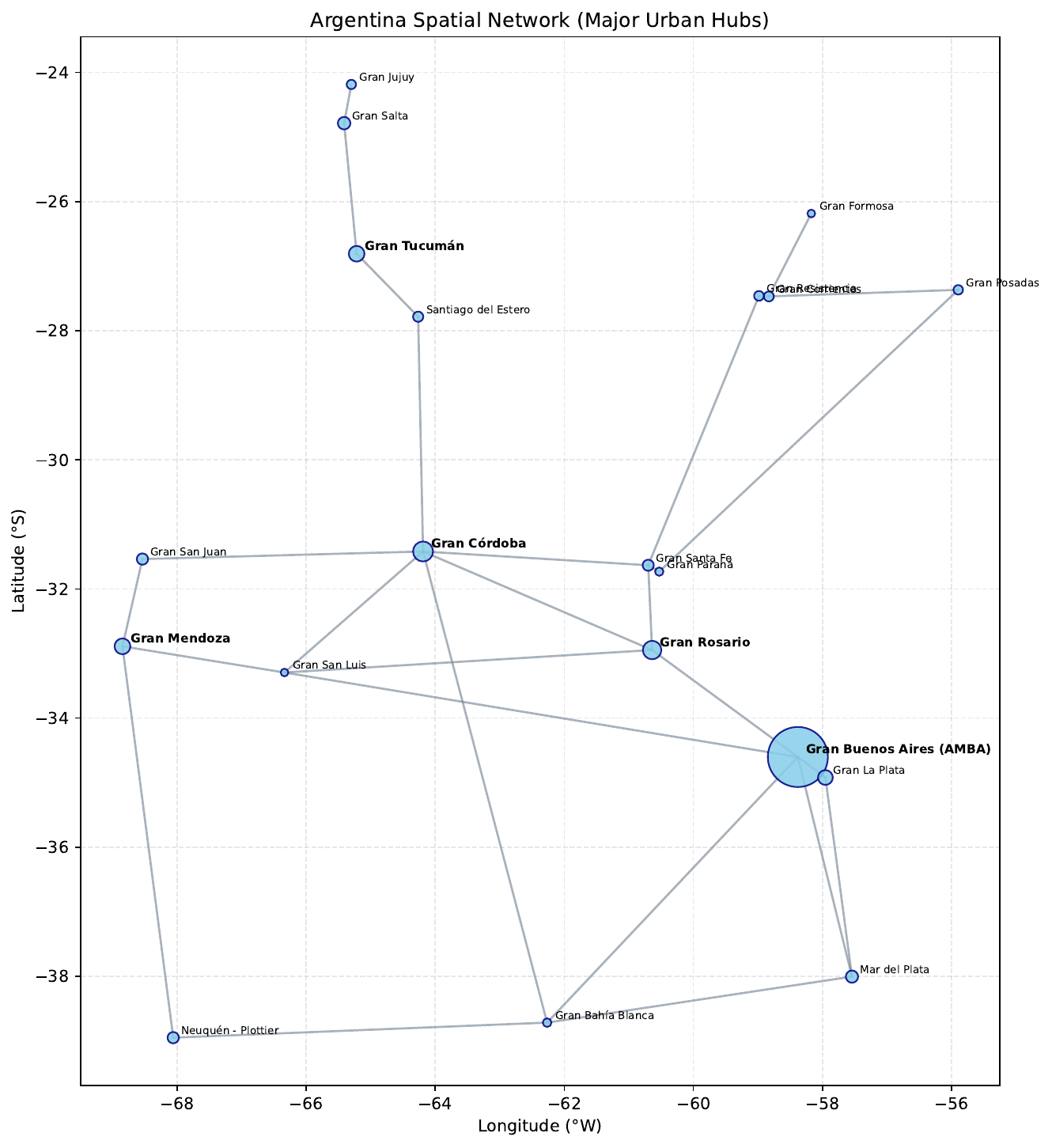}
\caption{Argentina spatial network of major urban areas and highway connections.}
\label{fig:arg_network}
\end{figure}

To build a pSIR model, we considered a spatial network representing the $20$ major urban areas across the country connected along major highway routes. The epidemic was assumed to originate in Greater Buenos Aires (AMBA), which was selected as the seed district. Figure~\ref{fig:arg_network} shows the resulting network, including the main urban areas and their connections.

Proceeding as in the previous examples, the districts were grouped into layers according to their graph distance from the seed district. The population distribution across these layers was then scaled to match the total population of the country at that time, $N = 45 \times 10^6$. This defines the layer population function $P(k)$, from which $\Pop(p_S)$ was computed using Eq.~\eqref{eq:Pop}.

Regarding the model parameters, we adopted $\gamma = 0.16$ and fitted $\beta$ and $\Delta t$ so that the simulated trajectory matched the observed data, noting that the remaining model parameters ($\Rinf$, $b_1$, $b_2$, etc.) follow directly from $\beta$ and $\gamma$. Assuming an underdetection factor of $1/11$, a nonlinear least-squares optimization (\texttt{scipy.optimize.curve\_fit} in Python) yielded $\Delta t = 14.99$ and $\beta = 0.2006$. Figure~\ref{fig:arg_fitted} compares the fitted trajectory with the reported data.

The results show that the aggregated pSIR model reproduces the overall trend of the observed epidemic curve, including its initial growth, peak, and subsequent decline. Overall, these results illustrate the ability of the aggregated pSIR model to capture the main features of real-world epidemic dynamics using a reduced representation of the underlying metapopulation network.

\FloatBarrier

\section{Conclusions}\label{sec:conclusions}

We introduced a novel compartmental model that approximates the aggregated dynamics of metapopulation SIR models. The model, called pSIR, consists of five compartments ($S$, $I$, $R$, $\pS$, and $\pR$), where $\pS$ contains the population living in districts where the epidemic has not yet arrived and $\pR$ contains the population living in districts where the epidemic has already finished.

The aggregated model contains six basic parameters: three corresponding to those of the classical SIR model ($\beta$, $\gamma$, and $N$) and three related to the network structure (number of districts $D$, mean connectivity $c$, and mean mobility $p$). In addition, the pSIR model incorporates a topology-dependent function $\Pop(p_S)$ that can be computed directly from the population size of each layer of the network.

We established fundamental dynamical and asymptotic properties of the pSIR trajectories and showed that, as in the standard SIR model, the basic reproduction number is given by $\Rz = \beta/\gamma$. Moreover, we proved that $\Rz$ coincides with the basic reproduction number of the underlying metapopulation model approximated by the pSIR system.

Three numerical examples were analyzed: two synthetic metapopulation models---first on a regular rectangular grid and then on a randomly generated spatial network with heterogeneous district populations---and a real-world application reproducing reported COVID-19 cases during the 2020 wave in Argentina using a network of the most populated urban areas and their interconnecting highways. In the real-data example, adjusting only two unknown parameters sufficed to fit the simulated trajectories to the observed incidence curve.

This work opens several directions for future research. The framework introduced here for the SIR model could be extended to other compartmental structures, such as SEIR or vector-borne SIR--SI models; while an SEIR extension appears straightforward, the SIR--SI case is more involved because both host and vector populations evolve in a spatially coordinated fashion. Another promising application is studying the effect of control interventions that restrict mobility between districts, where a low-dimensional model like pSIR can substantially simplify the design of feedback controllers and optimal control strategies compared with full metapopulation models.

\section*{Statements and Declarations}

\noindent\textbf{Data and Code Availability:} The Python source code used to generate the simulations in this study is publicly available at \url{https://github.com/CIFASIS/pSIR/}. The COVID-19 epidemiological data for Argentina are available at \url{https://covidstats.com.ar} and \url{https://fceia.unr.edu.ar/~kofman/files/nacion.csv}.

\medskip
\noindent\textbf{Competing Interests:} The authors declare that they have no competing interests.

\appendix

\section{Closed-Form Expression for $\Pop(p_S)$}\label{sec:closed_form}

Consider a regular rectangular grid of $M_1 \times M_2$ districts (with $M_1 \geq M_2$) where each interior district is connected to its four neighbors. When the $0$-th layer is a district located in a corner, $P(k)$ takes the form of Eq.~\eqref{eq:Pk}.

While $P(k) = N_d(k + 1)$, the cumulative layer population is
\begin{equation*}
P_c(k) = \sum_{j=0}^{k} N_d(j + 1) = \frac{N_d}{2}(k + 1)(k + 2).
\end{equation*}
Then, the midpoints $m_k = (P_c(k-1) + P_c(k))/2$ are given by
\begin{equation*}
m_k = P_c(k) - \frac{1}{2}P(k) = \frac{N_d}{2}(k+1)(k+2) - \frac{N_d}{2}(k+1) = \frac{N_d}{2}(k+1)^2,
\end{equation*}
which corresponds to the layer population $P(k) = N_d(k+1)$. From Eq.~\eqref{eq:Popdisc}, it follows that
\begin{equation*}
\Pop(N - m_k) = \Pop\!\left(N - \frac{N_d}{2}(k+1)^2\right) = P(k) = N_d(k + 1).
\end{equation*}
Setting $p_S = N - \frac{N_d}{2}(k+1)^2$, we obtain $N_d(k+1) = \sqrt{2N_d(N - \pS)}$, so that
\begin{equation*}
\Pop(\pS) = \sqrt{2N_d(N - \pS)}.
\end{equation*}
When $k$ is such that $P(k) = N_d M_2$, we directly obtain $\Pop(\pS) = N_d M_2$. For larger values of $k$, when $P(k) = N_d(M_1 + M_2 - 1 - k)$, a symmetric situation to that of $P(k) = N_d(k+1)$ occurs, yielding
\begin{equation*}
\Pop(\pS) = \sqrt{2N_d\pS}.
\end{equation*}
Combining the three regimes leads to Eq.~\eqref{eq:PopRect}. Although this expression was deduced at the discrete points $\pS = N - m_k$, it extends continuously to arbitrary values $\pS \in [0, N]$.

\section{Computing $\mathcal{R}_0$ for the Aggregated Model} \label{AppR0}

We follow the next-generation matrix method proposed in \cite{VanDenDriessche2002} to compute the basic reproduction number $\mathcal{R}_0$ for the aggregated pSIR model. To this end, we reorder the state variables in~\eqref{eq:pSIR} so that the infectious compartment $I$ appears first:
\begin{equation}\label{eq:reordered}
\begin{aligned}
\dot{I}(t) &= \beta \frac{S(t)I(t)}{S(t)+I(t)+R(t)} - \gamma I(t), \\
\dot{S}(t) &= g_S(I(t), p_S(t)) - \beta \frac{S(t)I(t)}{S(t)+I(t)+R(t)} - S_\infty g_R(S(t), I(t), R(t)), \\
\dot{R}(t) &= \gamma I(t) - R_\infty g_R(S(t), I(t), R(t)), \\
\dot{p}_S(t) &= -g_S(I(t), p_S(t)), \\
\dot{p}_R(t) &= g_R(S(t), I(t), R(t)),
\end{aligned}
\end{equation}
where $p_I(t) = S(t) + I(t) + R(t) > 0$ for all $t \ge 0$, $R_\infty = 1 - S_\infty$, $g_S(I, p_S) \ge 0$ denotes the spatial activation rate (proportional to $(I/I_{\mathrm{peak}})^2 \Pop(p_S)$ with $\Pop(0)=0$), and $g_R(S, I, R) \ge 0$ denotes the deactivation rate (satisfying $g_R(S, I, R) = 0$ whenever $R \le b_2(S+I)$).

For system~\eqref{eq:reordered}, we denote the state vector by $x = (I, S, R, p_S, p_R)^T$. The disease-free equilibria (DFE) with $p_I > 0$ take the form
\begin{equation}\label{eq:dfe}
x_0 = (0, S_0, 0, p_{S0}, p_{R0})^T,
\end{equation}
with $S_0 > 0$, $p_{S0}, p_{R0} \ge 0$, and $S_0 + p_{S0} + p_{R0} = N$.

Following the notation in \cite{VanDenDriessche2002}, let $x = (x_1, \dots, x_n)^T$ with $x_i \ge 0$ for $i = 1, \dots, 5$, sorted so that the first $m = 1$ compartment corresponds to infected individuals, and define the set of disease-free states
\begin{equation*}
X_s \coloneq \{ x \ge 0 : x_i = 0, \ i = 1, \dots, m \}.
\end{equation*}
We write the dynamics of each compartment in the form
\begin{equation}\label{eq:fv_decomp}
\dot{x}_i = f_i(x) = \mathcal{F}_i(x) - \mathcal{V}_i(x),
\end{equation}
where $\mathcal{F}_i(x)$ is the rate of appearance of new infections in compartment $i$, and $\mathcal{V}_i = \mathcal{V}_i^- - \mathcal{V}_i^+$, with $\mathcal{V}_i^-$ representing the rate of transfer of individuals out of compartment $i$ and $\mathcal{V}_i^+$ the rate of transfer into compartment $i$. Writing $f, \mathcal{F}, \mathcal{V}^-, \mathcal{V}^+$ for the vector functions with components $f_i, \mathcal{F}_i, \mathcal{V}_i^-, \mathcal{V}_i^+$ ($i = 1, \dots, 5$), we have
\begin{equation} \label{eq:fv_vectors1}
\mathcal{F}(x) =
\begin{pmatrix}
\beta \dfrac{SI}{S+I+R} \\[6pt]
0 \\[4pt]
0 \\[4pt]
0 \\[4pt]
0
\end{pmatrix},
\qquad
\mathcal{V}^-(x) =
\begin{pmatrix}
\gamma I \\[6pt]
\beta \dfrac{SI}{S+I+R} + S_\infty g_R(S,I,R) \\[6pt]
R_\infty g_R(S,I,R) \\[6pt]
g_S(I, p_S) \\[4pt]
0
\end{pmatrix},
\end{equation}
\begin{equation} \label{eq:fv_vectors2}
\mathcal{V}^+(x) =
\begin{pmatrix}
0 \\[6pt]
g_S(I, p_S) \\[6pt]
\gamma I \\[6pt]
0 \\[4pt]
g_R(S,I,R)
\end{pmatrix},
\end{equation}
with $\mathcal{V} = \mathcal{V}^- - \mathcal{V}^+$ and $f = \mathcal{F} - \mathcal{V}$. We verify that these functions satisfy the hypotheses (A1)--(A5) of \cite{VanDenDriessche2002}:
\begin{itemize}
\item[(A1)] If $x \ge 0$, then $\mathcal{F}_i, \mathcal{V}_i^+, \mathcal{V}_i^- \ge 0$ for $i = 1, \dots, 5$.
\item[(A2)] For each $i = 1, \dots, 5$, if $x_i = 0$ then $\mathcal{V}_i^- = 0$. In particular, if $x \in X_s$ then $\mathcal{V}_i^- = 0$ for $i = 1$.
\item[(A3)] $\mathcal{F}_i = 0$ if $i > m$.
\item[(A4)] If $x \in X_s$, then $\mathcal{F}_i(x) = 0$ and $\mathcal{V}_i^+(x) = 0$ for $i = 1$.
\item[(A5)] If $\mathcal{F}(x)$ is set to zero, then all eigenvalues of $Df(x_0)$ have negative real parts.
\end{itemize}

Hypotheses (A1), (A3), and (A4) follow immediately from Eqs.~\eqref{eq:fv_vectors1}--\eqref{eq:fv_vectors2}. Hypothesis (A2) is also satisfied: for $i = 1$, setting $I = 0$ yields $\mathcal{V}_1^- = 0$; for $i = 2$, $g_R(0,I,R) = 0$; for $i = 3$, $g_R(S,I,0) = 0$; and for $i = 4$, $\Pop(0) = 0$ implies $g_S(I,0) = 0$.

To verify (A5), we examine the linearization at $x_0$. Observe that in an open neighborhood of $x_0$ (where $S_0 > 0$ and $R_0 = 0 < b_2 S_0$), we have $R < b_2(S + I)$, so $g_R(S, I, R) \equiv 0$ identically and all its first partial derivatives vanish at $x_0$. Likewise, since $g_S(I, p_S)$ is quadratic in $I$, all first partial derivatives of $g_S$ vanish at $I = 0$. Computing the Jacobian matrices at $x_0$ therefore gives
\begin{equation}\label{eq:DF}
\begin{split}
D\mathcal{F}(x_0)
&= \left.
\begin{bmatrix}
\beta \dfrac{S(S+R)}{(S+I+R)^2} & \beta \dfrac{I(I+R)}{(S+I+R)^2} & -\beta \dfrac{SI}{(S+I+R)^2} & 0 & 0 \\[8pt]
0 & 0 & 0 & 0 & 0 \\
0 & 0 & 0 & 0 & 0 \\
0 & 0 & 0 & 0 & 0 \\
0 & 0 & 0 & 0 & 0
\end{bmatrix}
\right|_{x=x_0} \\
&=
\begin{bmatrix}
\beta & 0 & 0 & 0 & 0 \\
0 & 0 & 0 & 0 & 0 \\
0 & 0 & 0 & 0 & 0 \\
0 & 0 & 0 & 0 & 0 \\
0 & 0 & 0 & 0 & 0
\end{bmatrix},
\end{split}
\end{equation}
and
\begin{equation}\label{eq:DV}
D\mathcal{V}(x_0) = \left.
\begin{bmatrix}
\gamma & 0 & 0 & 0 & 0 \\[6pt]
\beta \dfrac{S(S+R)}{(S+I+R)^2} & \beta \dfrac{I(I+R)}{(S+I+R)^2} & -\beta \dfrac{SI}{(S+I+R)^2} & 0 & 0 \\[8pt]
-\gamma & 0 & 0 & 0 & 0 \\
0 & 0 & 0 & 0 & 0 \\
0 & 0 & 0 & 0 & 0
\end{bmatrix}
\right|_{x=x_0}
=
\begin{bmatrix}
\gamma & 0 & 0 & 0 & 0 \\
\beta & 0 & 0 & 0 & 0 \\
-\gamma & 0 & 0 & 0 & 0 \\
0 & 0 & 0 & 0 & 0 \\
0 & 0 & 0 & 0 & 0
\end{bmatrix}.
\end{equation}
When $\mathcal{F}(x) \equiv 0$, the Jacobian $Df(x_0) = -D\mathcal{V}(x_0)$ is lower-triangular with the single nonzero eigenvalue $-\gamma < 0$ associated with the infected compartment (and zero eigenvalues along the continuum of disease-free equilibria, as in the classical closed SIR model without vital dynamics).

Now, setting
\begin{equation}\label{eq:FV_def}
F \coloneq \left[\frac{\partial \mathcal{F}_i}{\partial x_j}(x_0)\right], \qquad V \coloneq \left[\frac{\partial \mathcal{V}_i}{\partial x_j}(x_0)\right], \qquad \text{with } 1 \le i, j \le m,
\end{equation}
in view of \cite[Lemma~1]{VanDenDriessche2002}, $F = \beta$ is non-negative and $V = \gamma$ is non-singular. The basic reproduction number for system~\eqref{eq:reordered} is then given by
\begin{equation}\label{eq:R0_def}
\mathcal{R}_0 \coloneq \rho(F V^{-1}),
\end{equation}
where $\rho(A)$ denotes the spectral radius of a matrix $A$. We have
\begin{equation}\label{eq:FV_vals}
F \coloneq \frac{\partial \mathcal{F}_1}{\partial I}(x_0) = \beta, \qquad V \coloneq \frac{\partial \mathcal{V}_1}{\partial I}(x_0) = \gamma,
\end{equation}
which leads to
\begin{equation}\label{eq:R0_result}
\mathcal{R}_0 = \frac{\beta}{\gamma},
\end{equation}
coinciding with the basic reproduction number of the standard SIR model (see e.g.~\cite{brauer2012mathematical}).

\section{Computing $\mathcal{R}_0$ for the Metapopulation Model}\label{AppR0network}

We consider the metapopulation model~\eqref{eq:SIRi} consisting of $D$ districts of equal population size $N_d \coloneq N/D$. Defining the effective infectious population in district $i$ as
\begin{equation}\label{eq:Ihat}
\hat{I}_i(t) \coloneq (1 - c_i p) I_i(t) + \sum_{j \in C_i} p I_j(t) = \sum_{j=1}^{D} p_{ij} I_j(t),
\end{equation}
system~\eqref{eq:SIRi} can be rewritten for $i = 1, \dots, D$ as
\begin{equation}\label{eq:SIRi_compact}
\begin{aligned}
\dot{S}_i(t) &= -\beta \frac{S_i(t)\hat{I}_i(t)}{N_d}, \\
\dot{I}_i(t) &= \beta \frac{S_i(t)\hat{I}_i(t)}{N_d} - \gamma I_i(t), \\
\dot{R}_i(t) &= \gamma I_i(t),
\end{aligned}
\end{equation}
where $P = [p_{ij}]_{1 \le i,j \le D}$ is the mobility matrix with entries
\begin{equation}\label{eq:P_entries}
p_{ij} \coloneq
\begin{cases}
1 - c_i p, & \text{if } j = i, \\
p, & \text{if } j \in C_i, \\
0, & \text{otherwise.}
\end{cases}
\end{equation}
Since $p$ represents the daily time proportion spent in each neighboring district and $c_i = |C_i|$ is the number of neighbors of district $i$, we have $0 \le c_i p \le 1$ for all $i = 1, \dots, D$. Hence $P$ has non-negative entries and each row sums to
\begin{equation*}
\sum_{j=1}^{D} p_{ij} = (1 - c_i p) + \sum_{j \in C_i} p = 1 - c_i p + c_i p = 1,
\end{equation*}
so $P$ is row-stochastic.

\begin{proposition}\label{prop:R0_network}
The basic reproduction number $\hat{\mathcal{R}}_0$ for the metapopulation model~\eqref{eq:SIRi} is
\begin{equation*}
\hat{\mathcal{R}}_0 = \frac{\beta}{\gamma}.
\end{equation*}
\end{proposition}

\begin{proof}
Following the calculations for $\mathcal{R}_0$ in Appendix~\ref{AppR0}, we order the state vector with the $m = D$ infected compartments first, i.e., $x = (I_1, \dots, I_D, S_1, \dots, S_D, R_1, \dots, R_D)^T$. Writing $\dot{x} = \mathcal{F}(x) - (\mathcal{V}^-(x) - \mathcal{V}^+(x))$, we have
\begin{equation}\label{eq:network_FV_vectors}
\mathcal{F}(x) =
\begin{pmatrix}
\beta \dfrac{S_1}{N_d} \sum_{j=1}^{D} p_{1j} I_j \\
\vdots \\
\beta \dfrac{S_D}{N_d} \sum_{j=1}^{D} p_{Dj} I_j \\
0 \\
\vdots \\
0 \\
0 \\
\vdots \\
0
\end{pmatrix},
\qquad
\mathcal{V}^-(x) =
\begin{pmatrix}
\gamma I_1 \\
\vdots \\
\gamma I_D \\
\beta \dfrac{S_1}{N_d} \sum_{j=1}^{D} p_{1j} I_j \\
\vdots \\
\beta \dfrac{S_D}{N_d} \sum_{j=1}^{D} p_{Dj} I_j \\
0 \\
\vdots \\
0
\end{pmatrix},
\qquad
\mathcal{V}^+(x) =
\begin{pmatrix}
0 \\
\vdots \\
0 \\
0 \\
\vdots \\
0 \\
\gamma I_1 \\
\vdots \\
\gamma I_D
\end{pmatrix}.
\end{equation}
Computing the $D \times D$ Jacobian matrices $F = \left[\frac{\partial \mathcal{F}_i}{\partial I_j}(x_0)\right]$ and $V = \left[\frac{\partial \mathcal{V}_i}{\partial I_j}(x_0)\right]$ for $1 \le i,j \le D$ at a disease-free equilibrium $x_0$, we get
\begin{equation}\label{eq:network_FV_matrices}
F =
\begin{bmatrix}
\beta \dfrac{S_{1,0}}{N_d} p_{11} & \beta \dfrac{S_{1,0}}{N_d} p_{12} & \cdots & \beta \dfrac{S_{1,0}}{N_d} p_{1D} \\[8pt]
\beta \dfrac{S_{2,0}}{N_d} p_{21} & \beta \dfrac{S_{2,0}}{N_d} p_{22} & \cdots & \beta \dfrac{S_{2,0}}{N_d} p_{2D} \\[6pt]
\vdots & \vdots & \ddots & \vdots \\[6pt]
\beta \dfrac{S_{D,0}}{N_d} p_{D1} & \beta \dfrac{S_{D,0}}{N_d} p_{D2} & \cdots & \beta \dfrac{S_{D,0}}{N_d} p_{DD}
\end{bmatrix},
\qquad
V =
\begin{bmatrix}
\gamma & \cdots & 0 \\
\vdots & \ddots & \vdots \\
0 & \cdots & \gamma
\end{bmatrix}.
\end{equation}
At a completely susceptible disease-free equilibrium, we have $I_{i,0} = 0$ and $R_{i,0} = 0$, and thus $S_{i,0} = N_d$ for all $i = 1, \dots, D$. Therefore,
\begin{equation}\label{eq:F_betaP}
F = \beta P.
\end{equation}
Consequently, since $V^{-1} = \frac{1}{\gamma} I_D$,
\begin{equation}\label{eq:R0_network}
\hat{\mathcal{R}}_0 = \rho(F V^{-1}) = \frac{\beta}{\gamma} \rho(P) = \frac{\beta}{\gamma},
\end{equation}
where $\rho(P) = 1$ holds because $P$ is row-stochastic (non-negative entries and each row sums to $1$).
\end{proof}

\bibliographystyle{plain}
\bibliography{references}

\end{document}